\documentclass[a4paper, 10pt]{article}
\RequirePackage{url}
\usepackage{setspace}
\usepackage{fancyhdr}
\usepackage{bm}
\usepackage{ulem}
\usepackage{cancel}
\usepackage{latexsym}
\usepackage{amssymb,amsmath,amsfonts,
bbm,pifont,upgreek,bbold,accents}  
\usepackage[colorlinks=true]{hyperref}
\hypersetup{urlcolor=blue, citecolor=red}

\usepackage[dvips]{graphicx}
\newtheorem{theorem}{Theorem}[section]
\newtheorem{definition}{Definition}[section]
\newtheorem{proposition}{Proposition}[section]
\newtheorem{lemma}{Lemma}[section]

\newtheorem{example}{Example}[section]

\newtheorem{remark}{Remark}[section]

\expandafter\chardef\csname pre amssym.def
at\endcsname=\the\catcode`\@
\catcode`\@=11
\def\undefine#1{\let#1\undefined}
\def\newsymbol#1#2#3#4#5{\let\next@\relax
 \ifnum#2=\@ne\let\next@\msafam@\else
 \ifnum#2=\tw@\let\next@\msbfam@\fi\fi
 \mathchardef#1="#3\next@#4#5}
\def\mathhexbox@#1#2#3{\relax
 \ifmmode\mathpalette{}{\m@th\mathchar"#1#2#3}%
 \else\leavevmode\hbox{$\m@th\mathchar"#1#2#3$}\fi}
\def\hexnumber@#1{\ifcase#1 0\or 1\or 2\or 3\or 4\or 5\or 6\or 7\or
8\or
 9\or A\or B\or C\or D\or E\or F\fi}
\ifcase\@ptsize
 \font\tenmsb=msbm10
 \font\sevenmsb=msbm7
 \font\fivemsb=msbm5
\or
 \font\tenmsb=msbm10 scaled \magstephalf
 \font\sevenmsb=msbm7 scaled \magstephalf
 \font\fivemsb=msbm5  scaled \magstephalf
\or
 \font\tenmsb=msbm10 scaled \magstep1
 \font\sevenmsb=msbm7 scaled \magstep1
 \font\fivemsb=msbm5 scaled \magstep1
\fi
\newfam\msbfam
\textfont\msbfam=\tenmsb
\scriptfont\msbfam=\sevenmsb
\scriptscriptfont\msbfam=\fivemsb
\edef\msbfam@{\hexnumber@\msbfam}
\def\Bbb#1{\fam\msbfam\relax#1}
\def\widehat#1{\setboxz@h{$\m@th#1$}%
 \ifdim\wdz@>\tw@ em\mathaccent"0\msbfam@5B{#1}%
 \else\mathaccent"0362{#1}\fi}
\def\widetilde#1{\setboxz@h{$\m@th#1$}%
 \ifdim\wdz@>\tw@ em\mathaccent"0\msbfam@5D{#1}%
 \else\mathaccent"0365{#1}\fi}

\def\RIfM@{\relax\ifmmode}
\def\nonmatherr@#1{\errmessage{\string#1\space allowed only in math mode}}
\def\Bbb{\RIfM@\expandafter\Bbb@\else
 \expandafter\nonmatherr@\expandafter\Bbb\fi}
\def\Bbb@#1{{\Bbb@@{#1}}}
\def\Bbb@@#1{\fam\msbfam\relax#1}
\def\setboxz@h{\setbox\z@\hbox}
\def\wdz@{\wd\z@}
\catcode`\@=\csname pre amssym.def at\endcsname

\begin{document}

\title{\bf 
Perihelion librations in  the secular three--body problem}

\author{ 
Gabriella Pinzari\thanks{{{I wish to thank the Editor--in--Chief Alain Goriely for his interest; Ugo Locatelli for 
highlighting discussions};   Qinbo Chen, Jerome Daquin and Sara Di Ruzza for their interest and doubly {Qinbo Chen for pointing out a mistake in a previous formulation of Proposition~\ref{holomorphy}.} This work has been supported by  the European Research Council [Grant number 677793 Stable and Chaotic Motions in the Planetary Problem]. Figures~\ref{figure1},~\ref{figure2} and~\ref{figure3} have been produced using the software {\sc mathematica}\textsuperscript{\textregistered}.
}}
\vspace{-.2truecm}
\\{\footnotesize Dipartimento di Matematica T. Levi--Civita}
\vspace{-.2truecm}
\\{\footnotesize via Trieste, 63, 35131, Padova (Italy)}
\vspace{-.2truecm}
\\{\scriptsize gabriella.pinzari@math.unipd.it}
}\date{February,  2020}
\maketitle

\begin{abstract}\footnotesize{
A normal form theory for non--quasi--periodic systems is combined with the special properties of the partially averaged Newtonian potential pointed out in~\cite{pinzari19}  to prove, in the averaged, planar three--body problem, the existence of a plenty of motions 
 where, periodically, the perihelion of the inner body affords librations about one equilibrium position and its ellipse squeezes to a segment before reversing its direction and  again decreasing its eccentricity ({\it perihelion librations}). 
}

\end{abstract}

\maketitle

\tableofcontents

\renewcommand{\theequation}{\arabic{equation}}
\setcounter{equation}{0}

\newpage
\section{Introduction}\label{setup}

{%{\bf\large 1.1}  
This paper deals with certain motions of   three point masses undergoing  Newtonian attraction.  More precisely, we study the case of two light bodies orbiting  their common center of mass (a ``binary asteroid system'')  while interacting with a heavier mass (a ``planet''), whose position is external to their trajectories. As no Newtonian interaction can be neglected, there  is no reason to claim that the system undergoes the  ``Keplerian'' approximation,  successfully used in the so--called ``planetary'' model\footnote{{The planetary model is the dynamical system of $1+n$ ($n\ge 2$) point masses undegoing Newtonian attraction, in the case when one of the masses (``sun'') is much heavier than the remaining ones (``planets''), which are of comparable size.}}\cite{arnold63, fejoz04, laskarR95, pinzari-th09, chierchiaPi11b}.  We think to the case  that the time of mutual revolution of the lighter bodies is much shorter than the time scale of the motions of the massive body and that the ratio $\varepsilon$ between the semi--axis of the ellipse of the asteroids and the position ray of the heavier one keeps to be less than $\frac{1}2$, so that collisions between the asteroids and the planet are not possible (a body moving on a Keplerian ellipse does not go beyond twice the semi--major axis from its focus). We look at the system from a reference frame centred with one of the asteroids or with their center of mass, so as to deal with an effective two--particle system, given by the other asteroid and the planet.  We shall refer as ``asteroidal ellipse'' the 
 instantaneous  ellipse of this asteroid,  focused on the center of the reference frame. We are interested to its motions.  
To simplify the analysis a little bit, 
we introduce three  assumptions. The main one  is based on the belief that, as long as the difference between the time scales persists, the system is ``well represented'' by a 
a certain ``averaged'' problem, which we call call  {\it secular problem}. We remark that such average is meant with respect to the proper time of the asteroid, so it should not be confused with the homonymous procedure often studied in the literature (e. g.~\cite{fejoz2016}), where the Keplerian approximation is used for two particles about their common sun, and the average is done with respect to {\it both} their mean anomalies. %We stress that we completely neglect the legitimate issue of providing  bounds to   $\mu$ and $\kappa$ in order to make the averaging assumption consistent, as well as the closely related (and certainly interesting) problem of removing it. 
 \\
 The secular system is simpler than the original problem, as we loose information concerning the position of the asteroid.
In particular,  collisions between the lighter particles are not observable.
The degrees of freedom of the system are the motions of the eccentricity (or of the angular momentum) and of the pericenter direction of the asteroidal ellipse and  the motions of the massive body. 
We associate to such system a certain  ``limiting system'' which is similarly defined, but with the massive body being firm. From now on, we refer to such limiting problem as  ``unperturbed'', and to the full secular problem as ``perturbed''. The terminology is here used with abuse, as we do not assume that the massive body has slow velocity in the full problem. For the unperturbed system,  only movements of the asteroidal ellipse occur.   By~\cite{pinzari19} the unperturbed problem turns out to be integrable, in the sense that it possesses a complete family of independent and commuting first integrals. More importantly, it reveals a surprising property, which we named   {\it renormalizable integrability}. Such property (recalled for completeness in Section~\ref{The secular Euler Hamiltonian} below) offers a remarkable shortcut to the knowledge of  movements of the asteroidal ellipse. By~\cite{pinzari19a}, in the case that the  interacting particles are constrained on a plane (this is actually our second assumption), and $\varepsilon<\frac{1}2$, there are two stable equilibria such that the pericenter direction of the asteroidal ellipse in the unperturbed problem  affords small oscillations about them, while its angular momentum oscillates about zero, affording a periodic change of sign. Physically, this means that  the asteroidal ellipse  is highly  eccentric at any time and, moreover,  there are two times (``squeezing times'') in a period of oscillation of the pericenter when the eccentricity is equal to one. Namely, at those times, the ellipse becomes a segment. After the squeezing times, the eccentricity of the ellipse decreases while the  sense of the motion is reversed. We call ``perihelion librations'' such kind of motions.
 The question which here we address is whether perihelion librations do persist in the full perturbed problem, when the massive body moves. We are able to give a positive answer to this question under our third assumption, which  consists in taking  the total angular momentum of the system (which is preserved during the motion) equal to zero. Under this assumption, the symmetries of the Hamiltonian ensure that the equilibria persist also in the perturbed problem, even though the integrability is lost. 
On the other side, such assumption is a source of  difficulty, as the angular momenta of the two particles will simultaneously vanish, and hence collisions of the heavy body with the center of the system are to be controlled.
 We shall formulate our result below, after we have introduced some mathematical tool.}

\begin{figure}
\center{ \includegraphics[height=5.0cm, width=9.0cm]{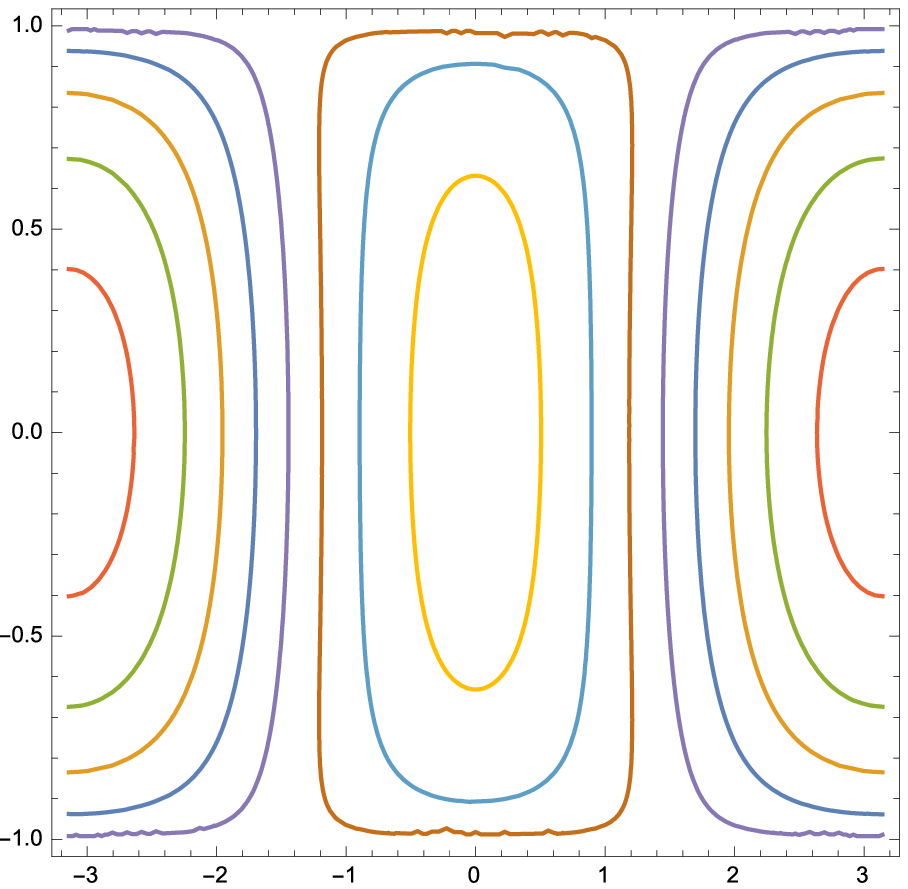}}
 \caption{Case $0<\varepsilon<\frac{1}2$.}\label{figure1}
 \center{\includegraphics[height=5.0cm, width=9.0cm]{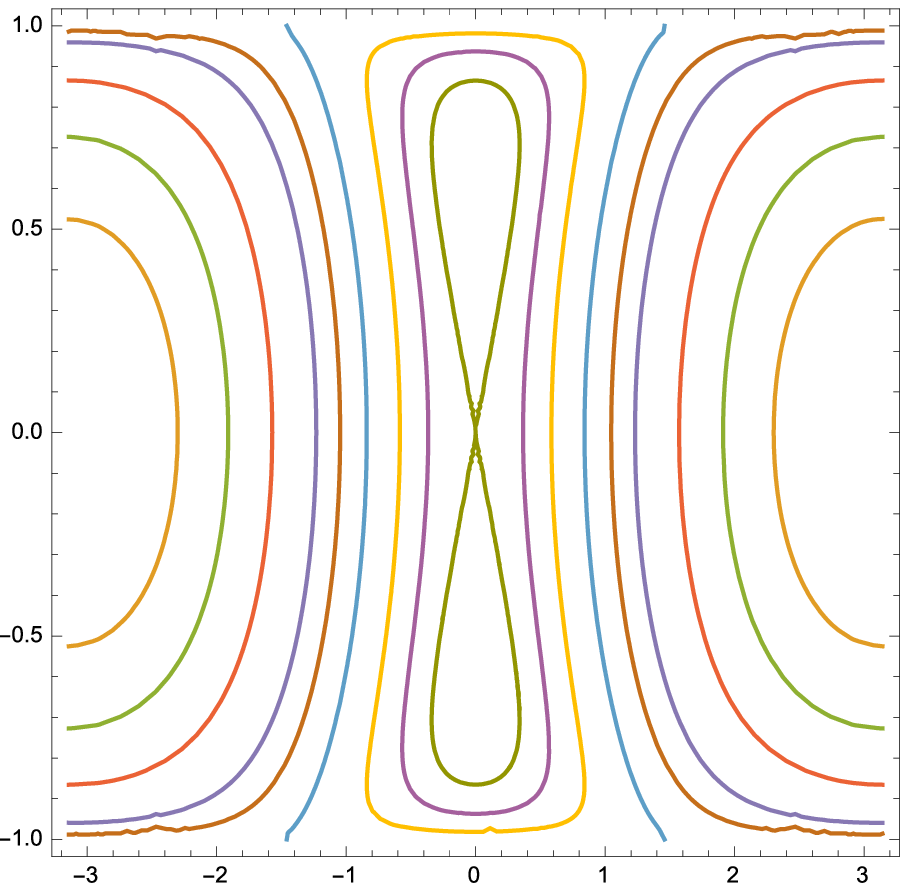}}  \caption{Case 
$\frac{1}2<\varepsilon<1$. }\label{figure2}
 \center{ \includegraphics[height=5.0cm, width=9.0cm]{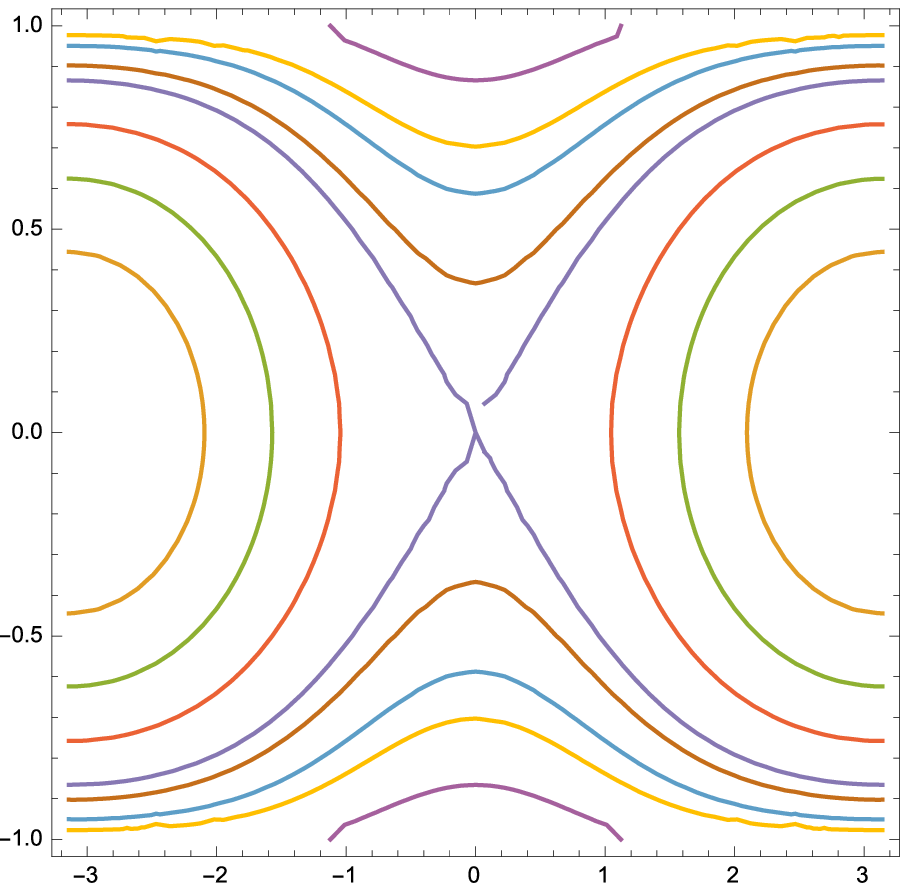}}  %
\caption{Case $\varepsilon>1$. }\label{figure3}
 \end{figure}

\vskip.1in
\noindent
{%\bf\large 1.2} {
In terms of Jacobi coordinates\footnote{{If (${\mathbf y}_0$, ${\mathbf x}_0$) (${\mathbf y}_1$, ${\mathbf x}_1$), (${\mathbf y}_2$,   ${\mathbf x}_2)$ are impulse--position coordinates of $m_0$, $\mu m_0$, $\kappa m_0$, respectively, by ``Jacobi coordinates'' one usually means  a linear, canonical  change (${\mathbf y}_0$, ${\mathbf y}_1$, ${\mathbf y}_2$, ${\mathbf x}_0$, ${\mathbf x}_1$, ${\mathbf x}_2)\in ({\mathbb R}^3)^6\to  ({\mathbf y}_{\rm c}$, ${\mathbf y}$, ${\mathbf y}'$, ${\mathbf x}_{\rm c}$, ${\mathbf x}$, ${\mathbf x}')\in ({\mathbb R}^3)^6$ defined so that  ${\mathbf x}_{\rm c}=({\mathbf x}_0+\mu{\mathbf x}_1+\kappa{\mathbf x}_2)(1+\mu+\kappa)^{-1}$ is the center of mass of the system,    while ${\mathbf x}$, ${\mathbf x}'$ are, respectively, the mutual position of ${\mathbf x}_1$ with respect to ${\mathbf x}_0$ and the position of ${\mathbf x}_2$ with respect to the center of mass of ${\mathbf x}_0$ and ${\mathbf x}_1$: ${\mathbf x}={\mathbf x}_1-{\mathbf x}_0$, ${\mathbf x}'={\mathbf x}_2-\frac{\mu}{1+\mu}{\mathbf x}_1-\frac{1}{1+\mu}{\mathbf x}_0$. The new impulses $({\mathbf y}_{\rm cm}$, ${\mathbf y}$, ${\mathbf y}'$) are uniquely defined by the  constraint of symplecticity. The new Hamiltonian turns to be ${\mathbf x}_{\rm c}$--independent, due to the conservation of ${\mathbf y}_{\rm c}$. The dependence on ${\mathbf y}_{\rm c}$ can be eliminated choosing (as it is always possible to do) a reference frame where ${\mathbf x}_{\rm c}\equiv{\mathbf 0}$. See e.g.,~\cite[\S 5.2--5.3]{giorgilli} for more details. }}    the three--body problem Hamiltonian with masses $m_0$, $\mu m_0$, $\kappa m_0$ is the translation--free function}
\begin{eqnarray*}
{{\rm H}_1}&=&{\frac{\|\mathbf y\|^2}{2m_0}\left(1+\frac{1}{\mu}\right)+\frac{\|\mathbf y'\|^2}{2m_0}\left(\frac{1}{1+\mu}+\frac{1}{\kappa}\right)-\frac{ \mu m^2_0}{\|\mathbf x\|}-\frac{\mu \kappa m_0^2}{\|\mathbf x'-\frac{1}{1+\mu}\mathbf x\|}-\frac{ \kappa m_0^2}{\|\mathbf x'+\frac{\mu}{1+\mu}\mathbf x\|}}
\end{eqnarray*}
{Here, $({\mathbf y}', {\mathbf y}, {\mathbf x}', {\mathbf x})\in ({\mathbb R}^3)^4$ (or $({\mathbb R}^2)^4$), $\|\cdot\|$ denotes Euclidean distance and the gravity constant has been taken equal to one, by a proper choice of the units system.}
{We rescale impulses and positions}
\begin{eqnarray} \label{resc}{\mathbf y\to \frac{\mu}{1+\mu}\mathbf y\ ,\quad \mathbf x\to {(1+\mu)}{\mathbf x}\ ,\quad \mathbf y'\to \mu\beta{\mathbf y'}\ ,\quad \mathbf x'\to \beta^{-1}\mathbf x'}\end{eqnarray} 
{multiply the Hamiltonian by $\frac{1+\mu}{\mu}$ (by a rescaling of time) and obtain
\begin{eqnarray} \label{Jac}{{\rm H}_1=\frac{\|\mathbf y\|^2}{2m_0}-\frac{m_0^2}{\|\mathbf x\|}+\gamma \left(\frac{\|\mathbf y'\|^2}{2m_0}
-\frac{\overline\beta}{\beta+\overline\beta}\frac{ m_0^2 }{\|\mathbf x'-\beta\mathbf x\|}
-\frac{\beta}{\beta+\overline\beta}\frac{ m_0^2 }{\|\mathbf x'+\overline\beta\mathbf x\|}
\right)}\end{eqnarray} 
 with}
\begin{eqnarray} \label{betagamma1}{\gamma=\frac{\kappa^3(1+\mu)^4}{\mu^3(1+\mu+\kappa)}\ ,\quad 
\beta=\frac{\kappa^2(1+\mu)^2}{\mu^2(1+\mu+\kappa)}\ ,\quad \overline\beta=\mu\beta}%\ ,\quad \overline\m:=\frac{1}{1+\mu}
\ .\end{eqnarray} 
{
Likewise, one might consider 
the problem written in the so--called $m_0$--centric\footnote{{If (${\mathbf y}_0$, ${\mathbf y}_1$, ${\mathbf y}_2$, ${\mathbf x}_0$, ${\mathbf x}_1$, ${\mathbf x}_2)$ are impulse--position coordinates of $m_0$, $\mu m_0$, $\kappa m_0$, respectively, the $m_0$--centric coordinates $({\mathbf Y}_0, {\mathbf y}, {\mathbf y}', {\mathbf X}_0, {\mathbf x}, {\mathbf x}')$ are defined via ${\mathbf X}_0={\mathbf x}_0$, ${\mathbf x}={\mathbf x}_1-{\mathbf x}_0$, ${\mathbf x}'={\mathbf x}_2-{\mathbf x}_0$ and by the symplecticity constraint.}} coordinates, and in this case the Hamiltonian is
}
\begin{eqnarray*}{{\rm H}_2=\frac{1}{2m_0}\left(1+\frac{1}{\mu}\right)\|{\mathbf y}\|^2+\frac{1}{2m_0}\left(1+\frac{1}{\kappa}\right)\|{\mathbf y}'\|^2-\frac{\mu m_0^2}{\|{\mathbf x}\|}-\frac{\kappa m_0^2}{\|{\mathbf x}'\|}-\frac{\mu \kappa m_0^2}{\|{\mathbf x}-{\mathbf x}'\|}+\frac{{\mathbf y}\cdot {\mathbf y}'}{m_0}\ .}\end{eqnarray*}
{We apply an analogue rescaling, but with 
}
\begin{eqnarray} \label{betagamma2}{\gamma=\frac{\kappa^3(1+\mu)^3}{\mu^3(1+\kappa)}\ ,\quad =\frac{\kappa^2(1+\mu)}{\mu^2(1+\kappa)}\ ,\quad \overline\beta=\mu\beta\ .}\end{eqnarray}  {We arrive at
\begin{eqnarray} \label{helio}{{\rm H}_2=\frac{\|{\mathbf y}\|^2}{2m_0}-\frac{m_0^2}{\|{\mathbf x}\|}+\gamma\left(
\frac{\|{\mathbf y}'\|^2}{2m_0}-\frac{\beta}{\beta+\overline\beta}\frac{ m_0^2}{\|{\mathbf x}'\|}-\frac{\overline\beta}{\beta+\overline\beta}\frac{m_0^2}{\|{\mathbf x}'-(\beta+\overline\beta){\mathbf x}\|}
\right)+\overline\beta \frac{{\mathbf y}'\cdot {\mathbf y}}{m_0}\ .}\end{eqnarray} 
{We remark that in the case, of our interest, that $\kappa\gg\mu\sim 1$,  the above definition of  Jacobi coordinates differs substantially from the usual one, because the barycentric reduction  begins with one of the two lighter masses, rather than with  the heavier one.  }
A similar observation holds for the  the $m_0$--centric reduction, which here is not centred on the most massive body, contrarily to the usual convention.
}
{We look at the Hamiltonians ${\rm H}_i$ in~(\ref{Jac}) and~(\ref{helio}).  As mentioned above, we make  three assumptions.}
 \begin{itemize}
\item[$A_1)$]  If $\ell$ is the mean anomaly associated to the Keplerian  motions of the term 
\begin{eqnarray} \label{Kep}\frac{\|\mathbf y\|^2}{2m_0}-\frac{m_0^2}{\|\mathbf x\|}=-\frac{m_0^5}{2{\Lambda}^2}\ ,\end{eqnarray} 
we replace the Hamiltonians~(\ref{Jac}) and~(\ref{helio}) with their respective  $\ell$--averages
\begin{eqnarray} \label{ovlH}\overline{\rm H}_i=-\frac{m_0^5}{2{\Lambda}^2}+\gamma\widehat{\rm H}_i\ ;\end{eqnarray} 
\item[$A_2)$] The coordinates ${\mathbf x}$, ${\mathbf x}'$ and the impulses ${\mathbf y}$, ${\mathbf y}'$ are constrained on the plane ${\mathbb R}^2$;
\item[$A_3)$] The total angular momentum $\mathbf C={\mathbf x}'\times {\mathbf y}'+{\mathbf x}\times {\mathbf y}$ of the system vanishes. \end{itemize}  
 {As mentioned above, the main assumption is $A_1)$. It allows us to exploit facts highlighted in~\cite{pinzari19, pinzari19a}, as now we describe.\\ Since the $\overline{\rm H}_i$'s in~(\ref{ovlH}) are $\ell$--independent, $\Lambda$ is a first integral, hence the 
 term $-\frac{m_0^5}{2{\Lambda}^2}$ may be neglected. After a further rescaling of time $t\to \gamma t$, we are led to look at the Hamiltonians $\widehat{\rm H}_i$ in~(\ref{ovlH}), which are given by   \begin{eqnarray} \label{secular}
\widehat{\rm H}_1&:=&\frac{\|{\mathbf y}'\|^2}{2 m_0}-\frac{m_0^2\overline\beta}{\beta+\overline\beta}{\rm U}%^+
_{\beta}-\frac{m_0^2\beta}{\beta+\overline\beta}{\rm U}%^-
_{-\overline\beta}\nonumber\\\nonumber\\ 
\widehat{\rm H}_2&:=&\frac{\|{\mathbf y}'\|^2}{2 m_0}-\frac{m_0^2\overline\beta}{\beta+\overline\beta}{\rm U}%^+
_{\beta+\overline\beta}-\frac{m_0^2\beta}{\beta+\overline\beta}\frac{1}{\|{\mathbf x}'\|}
\end{eqnarray} 
where  
\begin{eqnarray} \label{Usb}{\rm U}
_{\beta}:=\frac{1}{2\pi}\int_0^{2\pi}\frac{d\ell}{\|{\mathbf x}'-\beta{\mathbf x}(\ell)\|}\end{eqnarray} 
is the $\ell$--average of the Newtonian potential. Remark that 
  ${\mathbf y}(\ell)$ has  vanishing $\ell$--average\footnote{
For the unperturbed Keplerian motions, it is ${\mathbf y}={\rm m}\dot{\mathbf x}$, so
\begin{eqnarray*}\frac{1}{2\pi}\int_0^{2\pi}{\mathbf y}d\ell=\frac{1}{2\pi}\int_0^{2\pi}m_0\dot {\mathbf x}d\ell=\frac{1}{T}\int_0^{T}m_0 \dot {\mathbf x}dt=0\end{eqnarray*}
where $T=\frac{2\pi}{\omega_{\rm K}}$ and $\omega_{\rm K}=\frac{m_0^5}{{\Lambda}^3}$ is the Keplerian frequency.
}, so that the last term in~(\ref{helio}) does not survive. 
In the case of the planar problem, after the reduction of rotation invariance, the Hamiltonians $\widehat{\rm H}_i$ have two degrees of freedom. We use the following canonical coordinates
\begin{eqnarray*}
&&{\rm R}=\frac{{\mathbf y}'\cdot {\mathbf x}'}{\|{\mathbf y}'\|}\ ,\quad {\rm r}=\|{\mathbf x}'\|\nonumber\\
&& {\rm G}=\|{\mathbf x}\times {\mathbf y}\|\ ,\quad {\rm g}=\rm anomaly\ of\ {\mathbf P}\ with\ respect\ to\ a\ fixed\ direction\ {\mathbf a}
\end{eqnarray*}
where ${\mathbf P}$ is the perihelion of~(\ref{Kep}), and the direction ${\mathbf a}$, orthogonal to ${\mathbf x}\times{\mathbf y}$, will be specified later. Note the coordinates above are fit to describe motions of Keplerian elements for $({\mathbf y}, {\mathbf x})$, but not of
 $({\mathbf y}', {\mathbf x}')$. 
If ${\mathbf y}'$ is set to zero, the $\widehat{\rm H}_i$'s reduce to sums of averaged Newtonian potentials, which are integrable, as do not depend on ${\rm R}$.  The function
${\rm U}_\beta|_{\beta=1}$ has been thoroughly studied in~\cite{pinzari19}.
Its phase portrait in the plane $({\rm g}, {\rm G})$, while the ratio $\varepsilon=a/{\rm r}$ (where $a={\Lambda}^2/m_0^3$ is the semi--major axis associated to~(\ref{Kep})) varies, is as follows.
\begin{itemize}
\item[(i)] Case $0<\varepsilon<\frac{1}2$. There exist two centres, $(0, 0)$ and $(0, \pi)$ surrounded by librational motions (Figure~\ref{figure1}). 
\item[(ii)] Case $\frac{1}2<\varepsilon<1$. The equilibrium $(0, 0)$ becomes a saddle, with  its own separatrix (the light blue curve), while $(0, \pi)$ is still stable.  Two more equilibria appear on the ${\rm G}$--axis (Figure~\ref{figure2}). 
\item[(iii)] Case $\varepsilon>1$. The  equilibria on the ${\rm G}$--axis and the saddle  persist. There is the birth of rotational motions (Figure~\ref{figure3}).\end{itemize}
The purpose of this paper is to prove that motions of the kind (i) do persist in $\widehat{\rm H}_i$, when  ${\mathbf y}'\ne 0$. Note that we shall not require $\|{\mathbf y}'\|$ small. \\
To state the result, we  introduce the following quantities, which will be used as mass parameters, at the place of $\mu$ and $\kappa$: 
\begin{eqnarray} \label{b*}\beta_*:=\left\{\begin{array}{llll}
\displaystyle\frac{\beta\overline\beta}{\beta+\overline\beta}\qquad {\rm if}\quad i=1\\\\
\displaystyle\overline\beta\qquad\qquad {\rm if}\quad i=2
\end{array}\right.\qquad \beta^*:=\left\{\begin{array}{llll}\displaystyle \max\{\beta,  \overline\beta\}\quad {\rm if}\quad i=1\\\\\\
\displaystyle\beta+\overline\beta\qquad\qquad {\rm if}\quad i=2
\end{array}\right.\end{eqnarray}  
where $\beta$ and $\overline\beta$ are as in~(\ref{betagamma1}),~(\ref{betagamma2}), respectively. Observe that $\beta_*< \beta^*$ and  the case, of our interest,  $\kappa\gg\mu\sim 1$ corresponds to $\beta_*\sim \beta^*/2\gg 1$. \\  We shall prove the following result.}
{
\begin{theorem}[Perihelion librations about $(0,0)$]\label{main}
Fix an arbitrary neighbourhood ${\rm U}_0$ of $(0,0)$ and an arbitrary neighbourhood ${\rm V}_0$ of an unperturbed curve $\gamma_0(t)=({\rm G}_0(t), {\rm g}_0(t))\in {\rm U}_0$ in Figure~\ref{figure1}. Then it is possible to find six numbers $0<c<1$, $0<\beta_-<\beta_+$, $0<\alpha_-<\alpha_+$ $T>0$, such that, for any $\beta_-<\beta_*\le \beta^*<\beta_+$ 
the projections $\Gamma_0(t)=({\rm G}(t), {\rm g}(t))$ of
all the orbits $\Gamma(t)=({\rm R}(t), {\rm G}(t), {\rm r}(t), {\rm g}(t))$ of $\overline{\rm H}_1$, $\overline{\rm H}_2$ with initial datum $({\rm R}_0, {\rm r}_0, {\rm G}_0, {\rm g}_0)\in [\frac{1}{\sqrt{c\alpha_+}}, \frac{1}{\sqrt{c\alpha_-}}]\times [c\alpha_-, \alpha_+]\times {\rm U}_0$
belong to  ${\rm V}_0$ for all $0\le t\le T$. Moreover, the angle $\gamma(t)$  between the position ray of $\Gamma_0(t)$ and the ${\rm g}$--axis affords a variation larger than $2\pi$ during the time $T$.
\end{theorem}
}

{\smallskip\noindent}
{A similar statement concerning perihelion librations about $(0,\pi)$ holds true. The statement of Theorem~\ref{main}  deserves two remarks. The former regards  the motions involved, which are quasi--rectilinear, hence, close to be collisional. Generally speaking, for a three--body system composed of two asteroids and one planet,   three kinds of collisions are possible: (1) collisions between the two asteroids; (2) collision between one of the asteroids and the planet; (3) triple collision. The system under investigation is, as stressed above, an averaged problem derived from the full above problem. For this averaged problem collisions of kind (1) or (3) do not exist, as  the position of the asteroids is treated only in averaged meaning. 
 Collisions of the kind (2)  may exist, but they are to be intended as collisions between the planet and the average ellipse, rather than with a single particle on this ellipse.
 They are prevented by the assumption that the orbit of the planet is sufficiently far from the orbit of the asteroids. Namely, with a careful choice of the domain of the coordinates. Under this assumption, the averaged Hamiltonians $\widehat{\rm H}_i$ keep finite. Incidentally, their regularity is studied in Proposition~\ref{holomorphy}. During the proof of Theorem~\ref{main}, in Section 5, the trajectory of the massive planet is controlled to keep outside the trajectory of the asteroid for all the time of a perihelion libration.\\ The second remark concerns}  the thesis of Theorem~\ref{main}.   It holds in an {\it open} subset of phase space. In a sense, it recalls the statement  of Nekhorossev's Theorem~\cite{nehorosev77}. However, differently from it, Theorem~\ref{main} is not an application of   perturbation theory, nor it uses trapping arguments. The reason is the following. 
{In Section~\ref{Set out and analytic tools}, we shall see that the manifolds
\begin{eqnarray} \label{M0OLD}{\cal M}_0:=\{({\rm R}, {\rm G}, {\rm r}, {\rm g}):\ ({\rm G}, {\rm g})=(0, 0)\}\ ,\quad {\cal M}_\pi:=\{({\rm R}, {\rm G}, {\rm r}, {\rm g}):\ ({\rm G}, {\rm g})=(0, \pi)\}\end{eqnarray} 
are in fact invariant for $\widehat{\rm H}_i$. On such invariant manifold,  by $A_3)$, we have
\begin{eqnarray} \label{G=-GOLD}\|{\mathbf x}'\times {\mathbf y}'\|=\|-{\mathbf x}\times {\mathbf y}\|={\rm G}\ .\end{eqnarray} 
so
\begin{eqnarray*}\left.\|{\mathbf y}'\|^2\right|_{{\cal M}_0, {\cal M}_\pi}=\left.{\rm R}^2+\frac{{\rm G}^2}{{\rm r}^2}\right|_{{\cal M}_0, {\cal M}_\pi}= {\rm R}^2\ .\end{eqnarray*} Moreover, 
 the functions ${\rm U}_\beta$, ${\rm U}_{-\overline\beta}$, ${\rm U}_{\beta+\overline\beta}$ in~(\ref{secular}) are  asymptotic (as $\varepsilon\to 0$) to $ \frac{1}{{\rm r}}$.  
Hence, the motion of the coordinates $({\rm R}, {\rm r})$ on ${\cal M}_0$ and ${\cal M}_\pi$ is ruled by an Hamiltonian asymptotic  to
\begin{eqnarray*}\frac{{\rm R}^2}{2m_0}-\frac{m_0^2}{{\rm r}}\ .\end{eqnarray*}
This Hamiltonian generates unbounded (hence, non quasi--periodic) motions, for both positive and negative energies: For positive energies, both ${\rm R}$ and ${\rm r}$ are unbounded; for negative energies, only ${\rm R}$ is so.} In any case, these motions are not quasi--periodic and hence we cannot apply the machinery of perturbation theory. In Section~\ref{A normal form lemma without small divisors} we develop
 a theory suited to the case (see~\cite{fortunatiW16} for a result of the same kind). In this theory, no small denominators will arise, which is the reason why no trapping argument is needed.%
 %FINO QUI
 
 {\smallskip\noindent}
 Before switching to full statements and proofs, we quote three  open questions arising from the present setting.
\begin{itemize} \item[$\rm Q_1)$] Let us consider the cases $\frac{1}2<\varepsilon<1$ or $\varepsilon>1$ (Figures~\ref{figure2},~\ref{figure3}, respectively).  Does the separatrix   split  so as to produce chaotic dynamics in the partially averaged planar problem?
\item[$\rm Q_2)$] Again in  the cases above, let us consider the {\it full} three--body problem. It  has three degrees of freedom.  Does the separatrix   split  so as to produce Arnold instability~\cite{arnold64, Delshams2019}?
 \item[$\rm Q_3)$] 
 {What is the scenario  in the case of the spatial problem? }
 \end{itemize}
{This paper is organised as follows. In Section~\ref{review} we provide a review of the main results of~\cite{pinzari19, pinzari19a}. In particular, we recall the mathematical formulation of the mentioned renormalizable integrability and we carry from~\cite{pinzari19a} a set of action--angle like coordinates suited to our needs. In Section~\ref{A deeper look at the  planar case} we refine the analysis of ~\cite{pinzari19} to the case of the {\it planar} secular problem, in the region of phase space where $0<\varepsilon<\frac{1}2$. In this case we are able to obtain simpler formulae compared to ~\cite{pinzari19}
and hence  to study the regularity region of $\widehat{\rm H}_i$ completely. In Section~\ref{Set out and analytic tools} we state a normal form theory without small divisors (Theorem~\ref{NFL}), suited for a--periodic systems. The proof of Theorem~\ref{NFL} is deferred to Appendix~\ref{Time--one flows}. In Section~\ref{Completion of the proof of Theorem} we provide the proof of  Theorem~\ref{main}, as well as of a more precise version of it (Theorem~\ref{main1} below), as an application of Theorem~\ref{NFL}.}

\section{Review of the results of~\cite{pinzari19, pinzari19a}}\label{review}
\subsection{${\cal K}$ coordinates}\label{coordinates}

We describe  canonical
coordinates {suited to our problem}. 

{\smallskip\noindent}
We fix an arbitrary  orthonormal frame 
\begin{eqnarray*}{\rm F}_0:\qquad {\mathbf i}=\left(
\begin{array}{lll}
1\\
0\\
0
\end{array}
\right)\ ,\qquad {\mathbf j}=\left(
\begin{array}{lll}
0\\
1\\
0
\end{array}
\right)\ ,\qquad {\mathbf k}=\left(
\begin{array}{lll}
0\\
0\\
1
\end{array}
\right)\end{eqnarray*} in ${\mathbb R}^3$,
 that we call {\it inertial frame}

 {\smallskip\noindent}
 For given  $m_0>0$, 
 fix a  region of  phase space (i.e., a set of values of $({\mathbf y}', {\mathbf y},{\mathbf x}',  {\mathbf x})$) where the Kepler Hamiltonian~(\ref{Kep})
 takes negative values. Consider  the motion generated by~(\ref{Kep})  with initial datum $({\mathbf y}, {\mathbf x})$, and denote: 
 \begin{itemize}
 \item[{\tiny\textbullet}]
 $a$  the semi--major axis; 
 % \item[{\tiny\textbullet}]
 %$\ee$  the {\it eccentricity}; 
  \item[{\tiny\textbullet}] ${\mathbf P}$, with $\|{\mathbf P}\|=1$, the direction of perihelion, assuming the ellipse is not a circle;
    \item[{\tiny\textbullet}] $\ell$: the mean anomaly, defined, mod $2\pi$, as the area of the elliptic sector spanned by ${\mathbf x}$ from ${\mathbf P}$, normalized to $2\pi$.
%  \item[{\tiny\textbullet}]  the {\it true anomaly} $\n$, defined as $$\n=\arg(\cos\xi-\ee, \sqrt{1-\ee^2}\sin \xi)$$ with 
 %  \item[{\tiny\textbullet}]  the {\it eccentric anomaly} $\xi$, solving the {\it Kepler equation} $\xi-\ee\sin\xi=\ell$.
  %  \item[{\tiny\textbullet}] $\ee({\rm G})=\sqrt{1-\frac{{\rm G}^2}{{\Lambda}^2}}$ is the {\it eccentricity};
 %    \item[{\tiny\textbullet}]    the quantity $\varrho=1-\ee\cos\xi$ corresponding to the ratio $\frac{{\rm r}}{a}$.
    \end{itemize}

%{\smallskip\noindent} 
%Next, we introduce the following notations.
%\begin{itemize}
%     \item[{\tiny\textbullet}]  If ${\mathbf i}$, ${\mathbf k}\in {\mathbb R}^3$, with  ${\mathbf i}\perp {\mathbf k}$ and  ${\mathbf i}$, $\mathbf k$ not necessarily of length 1,  by ${\rm F}\sim({\mathbf i},\cdot,{\mathbf k})$, we mean  the orthonormal frame  ${\rm F}=(\frac{{\mathbf i}}{\|{\mathbf i}\|}, \frac{{\mathbf k}\times {\mathbf i}}{\|{\mathbf k}\times {\mathbf i}\|},\frac{{\mathbf k}}{\|{\mathbf k}\|})$.
%     \item[{\tiny\textbullet}] 
%Given a couple $({\rm F},{\rm F}')$ of orthonormal frames, where ${\rm F}\sim({\mathbf i},\cdot,{\mathbf k})$, ${\rm F}'\sim({\mathbf i}',\cdot,{\mathbf k}')$, with ${\mathbf k}\times{\mathbf k}'\ne {\mathbf 0}$, we write $${\rm F}\to^{(\YY, \XX,  x)}{\rm F}'$$ if ${\mathbf i}'={\mathbf k}\times{\mathbf k}'$ and
%$$\YY={\mathbf k}'\cdot \frac{\mathbf k}{\|{\mathbf k}\|}\ ,\qquad \XX=\|{\mathbf k}'\|\ ,\qquad  x=\a_{\mathbf k}({\mathbf i},{\mathbf i'})$$ 
%where $\a_{\mathbf k}({\mathbf i},{\mathbf i'})$ is the oriented angle  ${\mathbf i}$ to ${\mathbf i}'$, with respect to the counterclockwise orientation established by ${\mathbf k}$. 

%\end{itemize}

{\smallskip\noindent}
Denote also:
\begin{eqnarray*}{\mathbf M}={\mathbf x}\times {\mathbf y}\ ,\quad {\mathbf M}'={\mathbf x}'\times {\mathbf y}'\ ,\quad {\mathbf C}={\mathbf M}'+{\mathbf M}\ ,\end{eqnarray*}
where ``$\times$'' denotes skew--product in ${{\mathbb R}}^3$.
Observe the following relations
\begin{eqnarray} \label{orthogonality}{{\mathbf x}'}\cdot{{\mathbf C}}={{\mathbf x}'}\cdot{\big({\mathbf M}+{\mathbf M}'\big)}={{\mathbf x}'}\cdot{{\mathbf M}}\ ,\qquad {\mathbf P}\cdot {\mathbf M}=0\ ,\quad \|{\mathbf P}\|=1\ .\end{eqnarray} 
 Let
\begin{eqnarray} \label{nodes}{\mathbf i}_1:={\mathbf k}\times {\mathbf C}\ ,\quad {\mathbf i}_2:={\mathbf C}\times {\mathbf x}'\ ,\quad {\mathbf i}_3:={\mathbf x}'\times {\mathbf M}\ ,\quad {\mathbf i}_4:={\mathbf M}\times {\mathbf P}
\end{eqnarray} 
and  assume
\begin{eqnarray*}{\mathbf i}_{j}\ne 0\qquad j=1,\ 2,\ 3, 4\ .\end{eqnarray*}
Given three vectors ${\mathbf i}$, ${\mathbf i}'$ and ${\mathbf k}$, with ${\mathbf i}$, ${\mathbf i}'\perp{\mathbf k}$, we denote as $\alpha_{\mathbf k}({\mathbf i}, {\mathbf i}')$ the oriented angle  from ${\mathbf i}$ to ${\mathbf i}'$ relatively to the positive orientation established by ${\mathbf k}$.

{\smallskip\noindent}
We   define the coordinates \begin{eqnarray*}{{\cal K}}=(\Lambda, {\rm G}, \Theta, {\rm R}, {\rm C}, {\rm Z},  \ell, {\rm g}, \vartheta, \tt, \gamma, \zeta)\end{eqnarray*}
 as
\begin{eqnarray} \label{belline}
\left\{\begin{array}{llll}
 \displaystyle {\rm Z}={\mathbf C}\cdot {\rm k}\\ 
 \displaystyle {\rm C}=\|{\mathbf C}\|\\ 
  \displaystyle  {\rm R}=\frac{{\mathbf y}'\cdot {\mathbf x}'}{\|{\mathbf x}'\|}\\
\displaystyle  \Lambda=\sqrt{m_0^3 a}\\
 \displaystyle{\rm G}=\|{\mathbf M}\|\\
 \displaystyle \Theta=\frac{{\mathbf M}\cdot {\mathbf x}'}{\|{\mathbf x}'\|}\\
\end{array}\right.\qquad\qquad\qquad \left\{\begin{array}{llll}
 \displaystyle {\rm z}=\alpha_{{\rm k}}({\mathbf i}, {\mathbf i}_1)\phantom{{\rm Z}={\mathbf C}\cdot {\rm k}}\\ 
 \displaystyle \gamma=\alpha_{{\mathbf C}}({\mathbf i}_1, {\mathbf i}_2)\phantom{{\rm C}=\|{\mathbf C}}\\ 
  \displaystyle {\rm r}=\|{\mathbf x}'\|\phantom{{\rm R}=\frac{{\mathbf y}'\cdot {\mathbf x}'}{\|{\mathbf x}'\|}}\\ 
\displaystyle \ell=\textrm{\rm mean anomaly of } {\mathbf x}\ {\rm on\ \mathbb E}\phantom{\Lambda=\sqrt{m_0^3 a}}\\ 
\displaystyle {\rm g}=\alpha_{{\mathbf M}}({\mathbf i}_3, {\mathbf i}_4)\phantom{{\rm G}=\|{\mathbf M}}\\
  \displaystyle\vartheta=\alpha_{{\mathbf x}'}({\mathbf i}_2, {\mathbf i}_3)\phantom{\Theta=\frac{{\mathbf M}\cdot {\mathbf x}'}{\|{\mathbf x}'\|}}
\end{array}\right.
\end{eqnarray}

 % fino qui

{\smallskip\noindent}
The canonical character of ${{\cal K}}$ has been discussed in~\cite{pinzari19}, based on~\cite{pinzari13}.
 
{\smallskip\noindent}
 Using the formulae in the previous section, we provide the expressions of the following functions
\begin{eqnarray*}{\rm U}=\frac{1}{2\pi}\int_{0}^{2\pi}\frac{d\ell}{\|{\mathbf x}_{\cal K}'-{\mathbf x}_{\cal K}\|}\ ,\qquad {\rm E}=\|{\mathbf x}_{\cal K}\times {\mathbf y}_{\cal K}\|^2-m_0^3 {\rm e}_{\cal K} {\mathbf P}_{\cal K}\cdot {\mathbf x}_{\cal K}'\end{eqnarray*}
 (where  ${\mathbf x}_{\cal K}:={\mathbf x}\circ{\cal K}$, etc.) which will be mentioned in the next section. They are:

\begin{eqnarray} \label{EEJJ}
{\rm U}( \Lambda, {\rm G}, \Theta, {\rm r} , \ell, {\rm g})&=&\frac{1}{2\pi}\int_0^{2\pi}\frac{d\ell}{\sqrt{{{\rm r} }^2+2{\rm r}  a(\Lambda)\sqrt{1-\frac{\Theta^2}{{\rm G}^2}} {\rm p}(\Lambda, {\rm G}, \ell, {\rm g})+{a(\Lambda)}^2\varrho(\Lambda, {\rm G}, {\rm g})^2}}\nonumber\\
{\rm E}(\Lambda, {\rm G}, \Theta, {\rm r} , \ell, {\rm g})&=&{\rm G}^2+m_0^3{\rm r} \sqrt{1-\frac{\Theta^2}{{\rm G}^2}}\sqrt{1-\frac{{\rm G}^2}{{\Lambda}^2}}\cos{\rm g}
\end{eqnarray} 
where
 $a=a(\Lambda)$ the {\it semi--major axis}; ${\rm e}={\rm e}(\Lambda, {\rm G})$, the {\it eccentricity} of the ellipse,
$\varrho=\varrho(\Lambda, {\rm G}, \ell)$, ${\rm p}={\rm p}( \Lambda, {\rm G}, \ell, {\rm g})$ are defined as
\begin{eqnarray} \label{p}
a(\Lambda)&=&\frac{{\Lambda}^2}{m_0^3}\nonumber\\
{\rm e}(\Lambda, {\rm G})&:=&\sqrt{1-\frac{{\rm G}^2}{{\Lambda}^2}}\nonumber\\
 \varrho(\Lambda, {\rm G}, \ell)&:=&1-{\rm e}(\Lambda,  {\rm G})\cos\xi( \Lambda, {\rm G}, \ell)\nonumber\\
{\rm p}(\Lambda, {\rm G}, \ell, {\rm g})&:=&(\cos\xi( \Lambda, {\rm G}, \ell)-{\rm e}(\Lambda,  {\rm G}))\cos{\rm g}-\frac{{\rm G}}{\Lambda}\sin\xi( \Lambda, {\rm G}, \ell)\ ,\end{eqnarray} 
with $\xi=\xi(\Lambda, {\rm G}, \ell)$  the {\it eccentric anomaly}, defined as the solution of
{\it Kepler equation}
 \begin{eqnarray} \label{Kepler Equation}\xi-{\rm e}(\Lambda, {\rm G})\sin\xi=\ell\ .\end{eqnarray} 
Also these formulae have been discussed in~\cite{pinzari19a}.

\subsection{Renormalizable integrability}\label{The secular Euler Hamiltonian}

 We recall some results  concerning the functions ${\rm U}$ and ${\rm E}$ in~(\ref{EEJJ}). We refer to~\cite{pinzari19} for full details.

 \begin{definition}[{\cite[Definition 1]{pinzari19}}]\label{def: renorm integr}\rm Let $f$, $g$ be two  functions
of the form
\begin{eqnarray} \label{HJ}f(p, q, y, x)=\widehat f({\rm I}(p,q), y, x)\ ,\qquad g(p, q, y, x)=\widehat g({\rm I}(p,q), y, x)\end{eqnarray} 
where 
\begin{eqnarray} \label{D}(p, q, y, x)\in {\cal D}:={\cal B}\times U\end{eqnarray} 
with $ U\subset {\mathbb R}^2$, ${\cal B}\subset{\mathbb R}^{2n}$ open and connected, $(p,q)=$ $(p_1$, $\cdots$, $p_n$, $q_1$, $\cdots$, $q_n)$   {conjugate} coordinates with respect to the two--form $\o=dy\wedge dx+\sum_{i=1}^{n}dp_i\wedge dq_i$ and ${\rm I}(p,q)=({\rm I}_1(p,q), \cdots, {\rm I}_n(p,q))$, with
\begin{eqnarray*}{\rm I}_i:\ {\cal B}\to {\mathbb R}\ ,\qquad i=1,\cdots, n\end{eqnarray*}
pairwise Poisson commuting:
\begin{eqnarray} \label{comm}\big\{{\rm I}_i,{\rm I}_j\big\}=0\qquad \forall\ 1\le i<j\le n\qquad i=1,\cdots, n\ .\end{eqnarray} 
 We say that $f$ is {\it renormalizably integrable by $g$ via $\widetilde f$} (or {\it renormalizably integrable by $g$}, or simply {\it renormalizably integrable}),  if there exists a function  \begin{eqnarray*}\widetilde f:\qquad {\rm I}({\cal B})\times g(U)\to {\mathbb R}\ , \end{eqnarray*} such that
\begin{eqnarray} \label{renorm}f(p,q,y,x)=\widetilde f({\rm I}(p,q), \widehat g({\rm I}(p,q),y,x))\end{eqnarray} 
for all $(p, q, y, x)\in {\cal D}$.
\end{definition}
\begin{proposition}[{\cite[Proposition 4]{pinzari19}}]\label{rem}
If $f$  is renormalizably integrable by $g$, then:\item[{\rm (i)}] 
${\rm I}_1$, $\cdots$, ${\rm I}_n$  are first integrals to $f$ and $g$;
\item[{\rm (ii)}] 
$f$ and $g$ Poisson commute. \end{proposition}

{\smallskip\noindent}
Observe that,  if $f$ is renormalizably integrable via $g$, then, generically, their respective time laws for the coordinates $(y, x)$ are the same, up to rescaling the time. Formally:

\begin{proposition}[{\cite[Proposition 5]{pinzari19}}]\label{prop: fixed points}
{L}et $f$ be renormalizably integrable via $g$. Fix a value ${\rm I}_0$ for the integrals ${\rm I}$ and look at the motion of 
 $(y, x)$ under $f$ and $g$, on the manifold ${\rm I}={\rm I}_0$. For any fixed initial datum $(y_0, x_0)$, let 
 $g_0:=g({\rm I}_0, y_0, x_0)$. If
 $\o({\rm I}_0, g_0):=\partial_{g}\tilde f({\rm I}, g)|_{({\rm I}_0, g_0)}\ne 0$, the motion $(y^f(t), x^f(t))$ with initial datum$(y_0, x_0)$ under $f$ is related to the corresponding motion $(y^g(t), x^g(t))$  under  $g$ via
 \begin{eqnarray*} y^f(t)=y^g(\o({\rm I}_0, g_0) t)\ ,\quad x^f(t)=x^g(\o({\rm I}_0,  g_0) t)\end{eqnarray*}
In particular, under this condition, all the fixed points of $g$
in the plane $(y, x)$ are fixed point to $f$.
Values of $({\rm I}_0, g_0)$ for which $\o({\rm I}_0, g_0)= 0$ provide, in the plane $(y, x)$, curves of fixed points for $f$ (which are not necessarily curves of fixed points to $g$). 
\end{proposition}

{\smallskip\noindent}
We observe that ${\rm U}$ and ${\rm E}$
 have the form in~(\ref{HJ}), with ${\rm I}=({\rm I}_1, {\rm I}_2, {\rm I}_3)=({\rm r}, \Lambda, \Theta)$ verifying~(\ref{comm}) and $(y, x)=({\rm G}, {\rm g})$.

\begin{proposition}[{\cite[Proposition 6]{pinzari19}}]\label{main prop}
${\rm U}$ is renormalizably integrable via ${\rm E}$. Namely,
there exists a function ${\rm F}$ such that
\begin{eqnarray} \label{equality1}{\rm U}(\Lambda, {\rm G}, \Theta, {\rm G}, {\rm g}, {\rm r})={\rm F}\big(\Lambda, \Theta, {\rm r}, {\rm E}(\Lambda, {\rm G}, \Theta, {\rm G}, {\rm g}, {\rm r})\big)\ .\end{eqnarray} 
%The function ${\rm F}$ is given by
%\begin{eqnarray*}
%{\rm F}({\rm r}, \Lambda, \Theta, {\rm E})=\widetilde{\rm F}({\rm r}, a(\Lambda), {\cal E}({{\rm E}}), {\cI}(\Theta, {{\rm E}}))
%\end{eqnarray*}
%where
%\begin{eqnarray} \label{U} \widetilde{\rm F}({\rm r}, a,  \cE ,  \cI )=\frac{1}{2\pi}\int_{{\mathbb T}}\frac{(1- \cE \cos w)dw}{\sqrt{{\rm r}^2+a^2-2a({\rm r} \cI \sin w+a \cE \cos w)+a^2 \cE ^2\cos^2w}}\  ; \end{eqnarray} 
%\begin{eqnarray} \label{EI} {\cal E}({{\rm E}})=\frac{\sqrt{{\Lambdaambda}^2-{{\rm E}}}}{\Lambdaambda}\qquad {\cI}(\Theta, {{\rm E}})=\frac{\sqrt{{{\rm E}}-{\Theta^2}}}{\Lambda_2}\ .
%\end{eqnarray} 
\end{proposition}

%We  recall a result from \red{referenza}.
%Let
%\begin{eqnarray} \label{C}{\cal C}:\qquad(\ell,u,v)\in \cA\times{\mathbb T}\times V\to (\underline{\mathbf y}, \underline{\mathbf x})=({\mathbf y}', {\mathbf y}, {\mathbf x}', {\mathbf x})\in ({\mathbb R}^3)^4\end{eqnarray} 
% where $\cA$ is domain\footnote{{By ``domain'' we mean an open and connected set in ${\mathbb K}={\mathbb R}^m, \complex^m$.}} in ${\mathbb R}$,   $V$ is a domain in ${\mathbb R}^{10}$, ${\mathbb T}:={\mathbb R}/(2\pi{\mathbb Z})$, $(u, v)=\big((u_1, u_2, u_3, u_4, u_5)$, $(v_1, v_2, v_3, v_4, v_5)\big)$, 
%be a change of coordinates  which verifies
%$$d{\mathbf y}'\wedge d{\mathbf x}'+d{\mathbf y}\wedge d{\mathbf x}=d\Lambda\wedge d\ell+d u\wedge d v$$
%and
% \beq{2B}
%\left(\frac{\|{\mathbf y}\|^2}{2\mm}-\frac{\mm{\cal M}}{\|{\mathbf x}\|}\right)\circ{\cal C}=
%-\frac{m_0^5}{2{\Lambdaambda}^2}=:{\rm h}_{\rm Kep}(\Lambda)\ ,{\rm e}q 
%where $\mm$, ${\cal M}$ are fixed having assimed that the image of ${\cal C}$ in~\equ{C}
%is a domain of  $(\mathbf y, \mathbf x)$ where the left hand side of~\equ{2B} takes negative values.  Assume also that the  ellipse ${\mathbb E}$ generated by the two--body Hamiltonian~\equ{2B} has non--vanishing eccentricity. 
%We call any map of such kind {\it partial Kepler map}. 

{\smallskip\noindent}
%We  consider the ellipse generated by the ``Kepler Hamiltonian'' at left hand side in~\equ{2B} and denote as  ${\rm e}:=\sqrt{1-\frac{{\rm G}^2}{{\Lambdaambda}^2}}$ its eccentricity, where ${\rm G}:=\|{\mathbf M}\|$; ${\mathbf P}$, with $\|{\mathbf P}\|=1$ the direction of perihelion, assuming ${\rm e}\ne 0$. Recall that the eccentricity vector ${\mathbf L}$ in~\equ{CL} is related to ${\rm e}$ and ${\mathbf P}$ via  ${\mathbf L}=\mm^2{\cal M}'{\rm e}{\mathbf P}$, so that  the function $\EE$ in~\equ{kepler} becomes
%\begin{eqnarray} \label{cG}{{\rm E}}:={\rm G}^2-{\rm m}^2{\rm M} {\rm e}\,{\mathbf x}'\cdot {\mathbf P}\end{eqnarray} 

%{\smallskip\noindent}
% We consider the Taylor and  Taylor--Fourier expansions
%\begin{eqnarray*} {\rm U}(\Lambda, {\rm G}, \Theta, \ell, {\rm g}, {\rm r})=\sum_{n=0}^{+\infty} {\rm r}^{n} {\rm U}_{n}( \Theta, {\rm G}, {\rm g})=\sum_{n=0}^{+\infty}\sum_{m=0}^\infty {\rm r}^{n} {\rm U}_{nm}( \Theta, {\rm G})\cos m{\rm g}\ .\end{eqnarray*}
%In \red{referenza} we have discussed the following consequence of Theorem~\ref{partial integral}. 
% \begin{proposition}\label{HarrK}
% \begin{eqnarray} \label{Hprop1} {\rm U}_{nm}( \Theta, {\rm G})\equiv 0\qquad {\rm if}\qquad m\ge n-1\ ,\quad \forall \ n\ge 1\ .\end{eqnarray} 
% In particular, the term ${\rm U}_{1}$ vanishes identically and ${\rm U}_{2}$, called {the} dipolar term,  Poisson--commutes with  ${\rm G}$ (``Harrington property'').
% \end{proposition}
{\smallskip\noindent}
The phase portrait of ${\rm E}$ in the planar case is as in Figures~\ref{figure1},~\ref{figure2} and~\ref{figure3}, accordingly to the values of $\varepsilon$.
\subsection{Asymptotic  action--angle coordinates}\label{Asymptotic  action--angle coordinates}
In this section, we focus on the planar case, i.e., when
$\mathbf y'$, $\mathbf y$, $\mathbf x'$, $\mathbf x\in {\mathbb R}^2$. In that case,  the following  the 8--dimensional diffeomorphism replaces ${\cal K}$ in~(\ref{belline}):
\begin{eqnarray} \label{coord}
{\cal K}_0:\quad \left\{\begin{array}{llll}\displaystyle{\rm C}=\|\mathbf C\|\\
\displaystyle{\rm G}=\|\mathbf M\|\\
\displaystyle{\rm R}=\frac{\mathbf y'\cdot \mathbf x'}{\|\mathbf x'\|}\\
\displaystyle\Lambda={\rm m} \sqrt{{\cal M} a}
\end{array}\right.\qquad\qquad \left\{\begin{array}{llll}\displaystyle\gamma=\alpha_{\mathbf k}(\mathbf i, \mathbf x')+\frac{\pi}{2}\\
\displaystyle {\rm g}=\alpha_{\mathbf k}({\mathbf x'},\mathbf P)+\pi\\
\displaystyle {\rm r}=\|\mathbf x'\|\\
\displaystyle\ell={\rm mean\ anomaly\ of\ {\mathbf x}\ in\ \mathbb E}
\end{array}\right.
\end{eqnarray} 
${\cal K}_0$ may be regarded as the natural limit of ${\cal K}$, once $\Theta$ and $\vartheta$ are fixed to ($0$, $0$), ($0$, $\pi$) (which are the values they take in the planar case), respectively, and $({\rm Z}, {\rm z})$ are neglected.
The 
  functions ${\rm U}$    and ${\rm E}$ in~(\ref{EEJJ}) become
\begin{eqnarray} \label{UU**}{\rm U}(\Lambda, {\rm G}, {\rm g}, {\rm r})&:=&\frac{1}{2\pi}\int_0^{2\pi}\frac{d\ell}{\sqrt{{\rm r}^2+2  a(\Lambda) {\rm r}{\rm p}(\Lambda, {\rm G},\ell,  {\rm g})+ a(\Lambda)^2\varrho( \Lambda, {\rm G}, \ell)^2 }}\ ,\nonumber\\
{\rm E}(\Lambda, {\rm G}, {\rm g}, {\rm r})&=&{\rm G}^2+m_0^3{\rm r}\,\sqrt{1-\frac{{\rm G}^2}{{\Lambda}^2}}\cos{\rm g}\end{eqnarray} 
 where are $a(\Lambda)$, $\varrho( \Lambda, {\rm G}, \ell )$, ${\rm p}( \Lambda, {\rm G}, \ell, {\rm g})$ are as in~(\ref{p}). %We also let, for future need,
% \begin{eqnarray*}%{p}
%\varrho( \L, {\rm G}, \ell ):=\widehat\varrho( \L, {\rm G}, \xi( {\rm G}, \ell))\ ,\quad {\rm p}( {\rm G}, \ell, {\rm g}):=\widehat{\rm p}( {\rm G}, \xi( {\rm G}, \ell), {\rm g})
%\end{eqnarray*}
%where 
%\begin{eqnarray*}
% \widehat \varrho( {\rm G}, \xi):=1-{\rm e}(\red\L,  {\rm G})\cos\xi\ ,\quad 
%\widehat{\rm p}( {\rm G}, \xi, {\rm g}):=(\cos\xi-{\rm e}(\red\L,  {\rm G}))\cos{\rm g}-\frac{{\rm G}}{\Lambda}\sin\xi{\red{\sin\rm g}}
%\end{eqnarray*}
%and $\xi( {\rm G}, \ell)$, the ``eccentric anomaly'', solves the Kepler equation~\equ{Kepler equation}.

{\smallskip\noindent}
Unfortunately, the action--angle coordinates associated to ${\rm E}$ are not explicit, since they are defined via inversion of elliptic integrals. However, it is possible to define, explicitly, action--angle coordinates associated to the leading part of ${\rm E}$ in the case of 
of large ${\rm r}$:
\begin{eqnarray*}{\rm E}_1:= m_0^3{\rm r}\, \sqrt{1-\frac{{\rm G}^2}{{\Lambda}^2}} \cos{\rm g}\ .\end{eqnarray*}
As discussed in~\cite{pinzari19a}, these coordinates, denoted as\footnote{Beware not to confuse the coordinate $\gamma$ in~(\ref{Gg}) with its homonymous in~(\ref{belline}). The latter is a cyclic coordinate for the Hamiltonians ${\rm H}_i$ in~(\ref{Jac}) and~(\ref{helio}), and hence has no r\^ole in the paper.} $({\cal G}, \gamma)$, are defined via the canonical\footnote{
The change of coordinates discussed in~\cite{pinzari19a} is a little more general than~(\ref{Gg}), since it is a four--dimensional map composed by
(\ref{Gg}) and
\begin{eqnarray*}\left\{\begin{array}{llll}\displaystyle\Lambda={\cal L}\\\\
\displaystyle\ell=\l+\arg\left(\cos\gamma, \ \frac{{\cal L}}{|{\cal G}|}\sin\gamma\right)
\end{array}\right.\ .\end{eqnarray*}
In this paper we shall only use the two--dimensional projection~(\ref{Gg}).
} change

\begin{eqnarray} \label{Gg} 
\left\{\begin{array}{llll}
\displaystyle{\rm G}=\sqrt{{\Lambda}^2-{\cal G}^2}\cos\gamma\\
\displaystyle {\rm g}=-\tan^{-1}\left(\frac{\Lambda}{{\cal G}}\sqrt{1-\frac{{\cal G}^2}{{\Lambda}^2}}\sin\gamma\right)+k\pi\\ {\rm with}\quad k=\left\{\begin{array}{llll}0\ \ {\rm if}\ 0<{\cal G}<\Lambda\\
1\ \ {\rm if}\ -\Lambda<{\cal G}<0
\end{array}\right.%\ \sign\cos{\rm g}=\sign{\cal G}
\end{array}\right.
\end{eqnarray} 
for any fixed value of $\Lambda$. Observe that positive values of ${\cal G}$ (hence, $k=0$) provide coordinates with the image $({\rm G}, {\rm g})$ in a neighborhood of $(0,0)$; negative values  ($k=1$) are for $({\rm G}, {\rm g})$ in a neighborhood of $(0,\pi)$. 
Using these ``approximate'' coordinates,  one obtains the expression of ${\rm E}$ as a close--to--be--integrable system for large ${\rm r}$:
\begin{eqnarray} \label{EActionAngle}
{\rm E}=m_0^3{\rm r}\,\frac{{\cal G}}{\Lambda}+({\Lambda}^2-{\cal G}^2)\cos^2\gamma
\end{eqnarray} 
The coordinates $({\cal G}, \gamma)$ will be used in Section~\ref{Completion of the proof of Theorem}.
%\footnote{
%Here we have used the expansion
%\begin{eqnarray*}{\rm U}= -\frac{\mm{\cal M}'}{{\rm r}}-\mm{\cal M}'\frac{3}{2}\frac{a}{{\rm r}^2}\sqrt{1-\frac{{\rm G}^2}{{\Lambda}^2}}\cos{\rm g}+\OO({\rm r}^{-3})\ .\end{eqnarray*}
%Observe the proportionality between ${\rm E}_1$ and the term of order ${\rm r}^{-2}$ in the expansion of ${\rm U}$. Such proportionality persists also in the general case: 	\begin{eqnarray*}{\rm E}_1= m_0^3{\rm r}\,' \sqrt{1-\frac{{\rm G}^2}{{\Lambda}^2}} \sqrt{1-\frac{\Theta^2}{{\rm G}^2}}\cos{\rm g}\end{eqnarray*}
%and
%\begin{eqnarray*}{\rm U}= -\frac{\mm{\cal M}'}{{\rm r}}-\mm{\cal M}'\frac{3}{2}\frac{a}{{\rm r}^2}\sqrt{1-\frac{{\rm G}^2}{{\Lambda}^2}}\sqrt{1-\frac{\Theta^2}{{\rm G}^2}}\cos{\rm g}+\OO({\rm r}^{-3})\ .\end{eqnarray*}
%Relation~\equ{commutation} actually implies a weaker assertion: that the lower order terms in the respective expansions of ${\rm E}$ and ${\rm U}+\frac{\mm{\cal M}'}{{\rm r}}$ do commute. 
% }  of $\JJ_{\rm s}$ with respect to ${\rm r}$, which is
%\begin{eqnarray} \label{s, trunc}\JJ_{\rm s}({\rm r}, \L,  {\rm G},  {\rm g})=-\frac{m_0^5}{2{\Lambda}^2}-\frac{\mm{\cal M}'}{{\rm r}}-\mm{\cal M}'\frac{3}{2}\frac{a}{{\rm r}^2}\sqrt{1-\frac{{\rm G}^2}{{\Lambda}^2}}\cos{\rm g}+\OO({\rm r}^{-3})\ ,\end{eqnarray} 

 \section{A deeper look at the  planar case}\label{A deeper look at the  planar case}

In the planar case, the %function  ${\rm F}$ appearing in the 
relation~(\ref{equality1}) 
becomes very special.

{\smallskip\noindent}
Instead of ${\rm U}$ and ${\rm E}$, it is convenient to switch to the functions
\begin{eqnarray} \label{Ueps}&& 
\widehat{\rm U}_{\varepsilon}(\Lambda, {\rm G}, {\rm g}):=\frac{1}{2\pi}\int_0^{2\pi}\frac{d\ell}{\sqrt{1+2 \varepsilon{\rm p}(\Lambda, {\rm G},\ell,  {\rm g})+\varepsilon^2\varrho( \Lambda, {\rm G}, \ell)^2 }}\nonumber\\
&&\widehat{\rm E}_{\varepsilon}(\Lambda, {\rm G}, {\rm g}):=
\sqrt{1-\frac{{\rm G}^2}{{\Lambda}^2}}\cos{\rm g}+\varepsilon\frac{{\rm G}^2}{{\Lambda}^2}\ ,
\end{eqnarray} 
which are related to the previous ones via
\begin{eqnarray} \label{UUU}{\rm U}(\Lambda, {\rm G}, {\rm g}, {\rm r}):=\frac{1}{{\rm r}}\widehat{\rm U}_{\varepsilon(\Lambda, {\rm r})}(\Lambda, {\rm G}, {\rm g})\ ,\quad {\rm E}(\Lambda, {\rm G}, {\rm g}, {\rm r}):=m_0^3{\rm r}\,\widehat{\rm E}_{\varepsilon(\Lambda, {\rm r})}(\Lambda, {\rm G}, {\rm g})\end{eqnarray} 
with
\begin{eqnarray*}\varepsilon(\Lambda, {\rm r}):=%\frac{a}{{\rm r}}=
\frac{{\Lambda}^2}{m_0^3{\rm r}}=\frac{a(\Lambda)}{{\rm r}}\end{eqnarray*}
if $a=a(\Lambda)$ is as in~(\ref{belline}). 
We rewrite relation~(\ref{equality1}) as
\begin{eqnarray} \label{calU}\widehat{\rm U}_{\varepsilon(\Lambda, {\rm r})}(\Lambda, {\rm G}, {\rm g})=\widehat{\rm F}_{\varepsilon(\Lambda, {\rm r})}\big(\widehat{\rm E}_{\varepsilon(\Lambda, {\rm r})}(\Lambda, {\rm G}, {\rm g})\big)\end{eqnarray} 
Here we have used that $\widehat{\rm F}_\varepsilon$ does not depend explicitly on $\Lambda$, since both ${\rm U}$ and ${\rm E}$ depend on $\Lambda$ only via $\frac{{\rm G}}{\Lambda}$. 
We claim that
\begin{proposition}\label{prop: avU}
In the planar problem, {if $|\varepsilon|<\frac{1}2$},~(\ref{calU}) holds with
 \begin{eqnarray} \label{avU}
%{\rm F}({\rm r}, \L, t)=\widehat{\rm F}_{\varepsilon({\rm r})}(t)\quad{\rm with}\quad
\widehat{\rm F}_{\varepsilon}(t)=\frac{1}{2\pi}\int_0^{2\pi}\frac{(1-\cos\xi) d\xi}{\sqrt{1-2 \varepsilon  (1-\cos\xi) t+\varepsilon^2(1-\cos\xi)^2 }}\end{eqnarray} 
\end{proposition}

{\smallskip\noindent}
To prove Proposition~\ref{prop: avU}, we need to recall the following result from~\cite{pinzari19}.

 \begin{proposition}[{\cite[Theorem 4 and Remark 3]{pinzari19}}]\label{ren int} Let  $f({\rm P}, y, x)$ and $g({\rm P}, y, x)$  Poisson commute. Assume 
that, for any fixed ${\rm P}$, the level sets $\{(y, x):\ g({\rm P}, y, x)={\cal G}\}$ are  union of graphs
\begin{eqnarray*} y=g_{1, i}({\rm P}, x, {\cal G})\quad {\rm or} \quad x=g_{2, j}({\rm P}, y, {\cal G})\ .\end{eqnarray*}
Then $f$ renormalizably integrable by $g$ 
through $\widetilde f$, and $\widetilde f$ can be chosen to be
\begin{eqnarray} \label{fg}\widetilde f({\rm P},{\cal G})=f({\rm P}, y_j, g_{2, j}(y_j, {\cal G}))= f({\rm P}, g_{1, i}(x_i, {\cal G}), x_i)\ .\end{eqnarray} 
for some fixed $x_i$, $y_j$.
\end{proposition}

 {\smallskip\noindent}
{\bf Proof of Proposition~\ref{prop: avU}} We  apply Proposition~\ref{ren int} to the functions $\widehat{\rm U}_{\varepsilon(\Lambda, {\rm r})}(\Lambda, {\rm G}, {\rm g})$, $\widehat{\rm E}_{\varepsilon(\Lambda, {\rm r})}(\Lambda, {\rm G}, {\rm g})$ in~(\ref{UUU}). Such two functions  do commute since ${\rm U}(\Lambda, {\rm G}, {\rm g}, {\rm r})$ and ${\rm E}(\Lambda, {\rm G}, {\rm g}, {\rm r})$ do and they both commute with ${\rm r}$. Moreover, the level sets $\{(\Lambda, {\rm G}, {\rm g}):\ \widehat{\rm E}_{\varepsilon(\Lambda, {\rm r})}(\Lambda, {\rm G}, {\rm g})=t\}$ are graphs
\begin{eqnarray*}{\rm g}=\pm \cos^{-1}\left(\frac{t-\varepsilon(\Lambda, {\rm r})\frac{{\rm G}^2}{{\Lambda}^2}}{\sqrt{1-\frac{{\rm G}^2}{{\Lambda}^2}}}\right)+2j\pi=:g^\pm_j({\rm r}, \Lambda, {\rm G},  t)\ ,\quad j\in {\mathbb Z}\ .\end{eqnarray*} We use the formula in~(\ref{fg}) with ${\rm P}=(\Lambda, {\rm r})$, $g_{2, j}=g^\pm_j({\rm r}, \Lambda, {\rm G},  t)$, $f={\rm U}$ and $y_j={\rm G}_j=0$ for all $j$. When ${\rm G}=0$, the functions ${\rm g}^\pm_j$ take the value
\begin{eqnarray*}{\rm g}^\pm_j|_{{\rm G}=0}=\pm \cos^{-1}t+2j\pi\quad \forall\ {\rm r},\ \Lambda\end{eqnarray*} so, by~(\ref{fg}),
\begin{eqnarray*}
\widehat{\rm F}_{\varepsilon}(t)&=&\widehat{\rm U}_{\varepsilon}( 0, \pm\cos^{-1}t+2j\pi)\nonumber\\
&=&\frac{1}{2\pi}\int_0^{2\pi}\frac{(1-\cos\xi) d\xi}{\sqrt{1-2 \varepsilon  (1-\cos\xi) t+\varepsilon^2(1-\cos\xi)^2 }}\ .\quad \square
\end{eqnarray*}

{\smallskip\noindent}
 A
first consequence of the formula in~(\ref{avU})  is underlined  in the following

\begin{remark}\rm Equation~(\ref{avU}) can also be used to provide an expansion of $\widehat{\rm U}_\varepsilon$ about the equilibria $(0,0)$ and $(0, \pi)$. Indeed,  
$\widehat{\rm E}_{\varepsilon}(\Lambda, {\rm G}, {\rm g})$ takes the value $+1$ at  $(\Lambda, {\rm G}, {\rm g})=(0, 0)$, and the value $-1$ at $(\Lambda, {\rm G}, {\rm g})=(0, \pi)$. Therefore, an expansion of $\widehat{\rm F}_{\varepsilon}(t)$ about $\pm 1$ provides, via~(\ref{calU}), an expansion $\widehat{\rm U}_{\varepsilon}$ about the corresponding equilibrium. On the other hand, the value  of  $\widehat{\rm F}_{\varepsilon}(t)$ and of its derivatives at $t=+ 1$ or $t=- 1$ can be explicitly computed, using the residuus theorem.
For example, for $0<\varepsilon<\frac{1}2$, 
\begin{eqnarray} \label{Ieps}\widehat{\rm F}_{\varepsilon}(1)=\frac{1}{2\pi}\int_0^{2\pi}\frac{(1-\cos\xi)d\xi}{1-\varepsilon(1-\cos\xi)}=\frac{2}{\sqrt{1-2\varepsilon}(1+\sqrt{1-2\varepsilon})}\ .\end{eqnarray} 
\end{remark}

{\smallskip\noindent}
Another consequence of~(\ref{avU}) is:
%The equality $\cI_{-\varepsilon({\rm r})}(-t)=\widehat{\rm F}_{\varepsilon({\rm r})}(t)$
%also implies
%$$\widehat{\rm U}_{-\varepsilon({\rm r})}(\L, {\rm G}, {\rm g})=\cI_{-\varepsilon({\rm r})}(\widetilde\EE_{-\varepsilon({\rm r})}(\L, {\rm G}, {\rm g}))=\widehat{\rm F}_{\varepsilon({\rm r})}(-\widetilde\EE_{-\varepsilon({\rm r})}(\L, {\rm G}, {\rm g}))=\widehat{\rm F}_{\varepsilon({\rm r})}(\widetilde\EE_{\varepsilon({\rm r})}( {\rm G}, {\rm g}+\p))\ .$$

\begin{proposition}\label{holomorphy}
{Let $|\varepsilon|<\frac{1}2$.} The complex function $(\varepsilon, t)\to \widehat{\rm F}_{\varepsilon}(t)$  looses its holomorphy if and only if
\begin{eqnarray} \label{sing}
{4\varepsilon^2-4\varepsilon t+1=0}\ .
\end{eqnarray}

\end{proposition}

\par\medskip\noindent{\bf Proof\ } {
We equivalently write
 \begin{eqnarray*}
%{\rm F}({\rm r}, \L, t)=\widehat{\rm F}_{\varepsilon({\rm r})}(t)\quad{\rm with}\quad
\widehat{\rm F}_{\varepsilon}(t)=\frac{1}{\pi}\int_0^{\pi}\frac{(1-\cos\xi) d\xi}{\sqrt{1-2 \varepsilon  (1-\cos\xi) t+\varepsilon^2(1-\cos\xi)^2 }}\ .\end{eqnarray*}
Then we change variable, letting
$x=1-\cos\xi$. The integral becomes
 \begin{eqnarray*}
%{\rm F}({\rm r}, \L, t)=\widehat{\rm F}_{\varepsilon({\rm r})}(t)\quad{\rm with}\quad
\widehat{\rm F}_{\varepsilon}(t)=\frac{1}{\pi}\int_0^{2}\frac{xdx}{\sqrt{x(2-x)}\sqrt{1-2 \varepsilon  x t+\varepsilon^2x^2 }}\ .\end{eqnarray*}
The only possibility it diverges is that there are two coinciding roots of the denominator on the real interval $[0, 2]$. The polynomial under the second square root has not coinciding  roots on $[0, 2]$ for $|\varepsilon|<\frac{1}2$ and never has the root $x=0$. The only possibility is that it has the root $x=2$. This happens when $t=\varepsilon+\frac{1}{4\varepsilon}$. $\qquad \square$}

	 %$|\varepsilon|<\frac{1}2$ and $t\in \complex$. 
{\smallskip\noindent} 
%\end{corollary}
%Corollary~\ref{cor: holomorphy} reflects on ${\rm U}$  as follows.

%\begin{proposition}\label{holomorphyU+}   The function ${\rm U}(\L, {\rm G}, {\rm g}, {\rm r})$  in~\equ{UU**} looses its holomorphy either when
%\begin{eqnarray} \label{no holomorphy}
%&&\left(\frac{{\Lambda}^2}{m_0^3{\rm r}}, \frac{\EE(\L, {\rm G}, {\rm g}, {\rm r})}{m_0^3{\rm r}}\right)\in {\mathbb R}^2\ ,\quad |\EE(\L, {\rm G}, {\rm g}, {\rm r})|\gammae m_0^3|{\rm r}|\ ,\quad  \sign \frac{{\Lambda}^2\EE(\L, {\rm G}, {\rm g}, {\rm r})}{m_0^6{\rm r}^2}=+1\nonumber\\
%&& \frac{1}2\left|\frac{\EE(\L, {\rm G}, {\rm g}, {\rm r})}{{\Lambda}^2}\right|\le 1\end{eqnarray} 
%or when $\EE(\L, {\rm G}, {\rm g}, {\rm r})$ does, namely, when ${\rm G}=\L$.\end{proposition}
 \begin{remark}\rm
The formula in~(\ref{avU}) (and its consequences below) is pretty specific for the planar case.  In~\cite[Equation (49)]{pinzari19} we proposed a general formula (holding for the planar and the spatial case),  which, unfortunately, does not seem equally exploitable. 
 \end{remark}

\section{Set out and analytic tools}\label{Set out and analytic tools}
For definiteness, we refer to perihelion librations about $(0, 0)$; the case $(0, \pi)$ being specular.

{\smallskip\noindent} 
Using  the identity in~(\ref{calU}), 
the assumption  $A_3$ in the introduction,  which gives \begin{eqnarray} \label{G}\|{\mathbf x}'_{{\cal K}_0}\times {\mathbf y}'_{{\cal K}_0}\|=\|{\mathbf x}_{{\cal K}_0}\times {\mathbf y}_{{\cal K}_0}\|={\rm G}\ .\end{eqnarray} 
and
 the relation
 \begin{eqnarray} \label{G1}
\|{\mathbf y}'_{{\cal K}_0}\|^2&=&\frac{|{\mathbf y}'_{{\cal K}_0}\cdot\mathbf x'_{{\cal K}_0}|^2}{\|\mathbf x'_{{\cal K}_0}\|^2}+\frac{\|{\mathbf y}'_{{\cal K}_0}\times\mathbf x'_{{\cal K}_0}\|^2}{\|\mathbf x'_{{\cal K}_0}\|^2}={\rm R}^2+\frac{{\rm G}^2}{{\rm r}^2}
\end{eqnarray} 
we rewrite the Hamiltonians~(\ref{secular}) as

\begin{eqnarray} \label{VVV}
\widehat{\rm H}_1({\rm R}, {\rm G}, {\rm r}, {\rm g})&=&\frac{{\rm R}^2}{2m_0}+\frac{{\rm G}^2}{2m_0{\rm r}^2}-\frac{\overline\beta}{\beta+\overline\beta}\frac{m_0^2}{{\rm r}}\widehat{\rm F}_{\beta\varepsilon({\rm r})}\Big(\widehat{\rm E}_{\beta\varepsilon({\rm r})}({\rm G}, {\rm g})\Big)\nonumber\\
&&-\frac{\beta}{\beta+\overline\beta}\frac{m_0^2}{{\rm r}}\widehat{\rm F}_{-\overline\beta\varepsilon({\rm r})}\Big(\widehat{\rm E}_{-\overline\beta\varepsilon({\rm r})}({\rm G}, {\rm g})\Big)\nonumber\\\nonumber\\
\widehat{\rm H}_2({\rm R}, {\rm G}, {\rm r}, {\rm g})&=&\frac{{\rm R}^2}{2m_0}+\frac{{\rm G}^2}{2m_0{\rm r}^2}-\frac{\overline\beta}{\beta+\overline\beta}\frac{m_0^2}{{\rm r}}\widehat{\rm F}_{(\beta+\overline\beta)\varepsilon({\rm r})}\Big(\widehat{\rm E}_{(\beta+\overline\beta)\varepsilon({\rm r})}({\rm G}, {\rm g})\Big)\nonumber\\
&&-\frac{\beta}{\beta+\overline\beta}\frac{m_0^2}{{\rm r}}
\end{eqnarray} 
 having neglected to write (as well as we shall do below) the dependence on $\Lambda$.

{\smallskip\noindent} 
We suddenly remark that $\widehat{\rm H}_1$ and $\widehat{\rm H}_2$ are both even with respect to ${\rm G}$ and ${\rm g}$ separately, because so are the functions $\widehat{\rm E}_\varepsilon({\rm G}, {\rm g})$ and the\footnote{Note that this circumstance would not hold without assuming $A_3$, since, in that case, instead of the term $\frac{{\rm G}^2}{{\rm r}^2}$, we would have 
$\frac{({\rm C}-{\rm G})^2}{{\rm r}^2}$ or $\frac{({\rm C}+{\rm G})^2}{{\rm r}^2}$ (depending on the verse of rotation of ${\mathbf x}'$),
which would break the symmetry.} term $\frac{{\rm G}^2}{2m_0{\rm r}^2}$.
Then the  manifolds {\begin{eqnarray} \label{M0}{\cal M}_0:=\{({\rm R}, {\rm G}, {\rm r}, {\rm g}):\ ({\rm G}, {\rm g})=(0, 0)\}\ ,\quad {\cal M}_\pi:=\{({\rm R}, {\rm G}, {\rm r}, {\rm g}):\ ({\rm G}, {\rm g})=(0, \pi)\}\end{eqnarray} }
which, by the discussions of Section~\ref{review}, are invariant for 
$\widehat{\rm E}_\varepsilon({\rm G}, {\rm g})$,
keep to be so also for $\widehat{\rm H}_i$. We focus on ${\cal M}_0$. On ${\cal M}_0$, the motions of the coordinates $({\rm R}, {\rm r})$ are governed by 
the Hamiltonians

%$$\widehat{\rm H}_i=\hh_i+f_i+k+g_i\ ,\qquad i=1,\ 2$$
%with
\begin{eqnarray} \label{hi}{\rm h}_i({\rm R}, {\rm r})=\frac{{\rm R}^2}{2m_0}+{\rm V}_i({\rm r})\end{eqnarray} 
where (using~(\ref{Ieps}) and dehomogenizating)
\begin{eqnarray} \label{Vi}
{\rm V}_1({\rm r})&=&%\frac{{\rm R}^2}{2m_0}
-\frac{\overline\beta}{\beta+\overline\beta}\frac{2m_0^2}{\sqrt{{\rm r}-2\b a}\left(\sqrt{\rm r}+\sqrt{{\rm r}-2\b a}\right)}-\frac{\beta}{\beta+\overline\beta}\frac{2 m_0^2}{\sqrt{{\rm r}+2\overline\beta a}\left(\sqrt{\rm r}+\sqrt{{\rm r}+2\overline\beta a}\right)}
\nonumber\\
{\rm V}_2({\rm r})&=&%\frac{{\rm R}^2}{2m_0}
-\frac{\overline\beta}{\beta+\overline\beta}\frac{2m_0^2}{\sqrt{{\rm r}-2(\beta+\overline\beta)a}\left(\sqrt{\rm r}+\sqrt{{\rm r}-2(\beta+\overline\beta)a}\right)}-\frac{\beta}{\beta+\overline\beta}\frac{ m_0^2}{{\rm r}}
\end{eqnarray} 
The ``potentials'' ${\rm V}_1$ and ${\rm V}_2$ are well defined and increasing from $-\infty$ to $0$ for  ${\rm r}>2\b a$, ${\rm r}>2(\beta+\overline\beta)$, respectively, so, action--angle coordinates do not exist.  In other words, closely to ${\cal M}_0$, $\widehat{\rm H}_i$ has not close to an integrable system in the sense of Liouville--Arnold~\cite{arnold63a} and hence  standard perturbation theory does not apply.
In the next section, we develop an analytic theory  suited to this case. 
It will be used to ``decouple'' the  the Hamiltonians. 
\subsection{A normal form theory without quasi--periodic unperturbed  motions}\label{A normal form lemma without small divisors}

In this section we describe a procedure for eliminating the angles\footnote{Note  that the procedure described in this section does not seem to be related to~\cite[\S 6.4.4]{arnoldKS06}, for the lack of slow--fast couples.} $\bm\varphi$ at high orders, given Hamiltonian of the form
 \begin{eqnarray} \label{h0f0}{\rm H}({\mathbf I}, \bm\varphi, {\mathbf p}, {\mathbf q}, y, x)={\rm h}({\mathbf I},{\mathbf J}({\mathbf p}, {\mathbf q}), y)+f({\mathbf I}, \bm\varphi, {\mathbf p}, {\mathbf q}, y, x)\end{eqnarray} 
%$${\rm h}(y,{\mathbf I},{\mathbf J})={\mathbf J}({\mathbf I}, y)+{\mathbf P}(y,{\mathbf I},{\mathbf J})\ ,\qquad {\mathbf P}(y,{\mathbf I},0)\equiv0$$
which we assume to be holomorphic on the neighbourhood
\begin{eqnarray*}{\mathbb P}_{\rho, s, \delta, r, \xi}={\mathbb I}_\rho  \times {{\mathbb T}}^n_s\times {\mathbb B}_{\delta}\times{\mathbb Y}_r\times  {\mathbb X}_\xi\supset {\mathbb P}={\mathbb I}  \times {{\mathbb T}}^n\times {\mathbb B}\times{\mathbb Y}\times  {\mathbb X}\ ,\end{eqnarray*}
and
\begin{eqnarray*}{\mathbf J}({\mathbf p}, {\mathbf q})=(p_1q_1,\cdots, p_mq_m)\ .\end{eqnarray*}
Here, ${\mathbb I}\subset {\mathbb R}^n$, ${\mathbb B}\subset {\mathbb R}^{2m}$, ${\mathbb Y}\subset {\mathbb R}$, ${\mathbb X}\subset {\mathbb R}$ are open and connected; $\mathbb T=\mathbb R/(2\pi\mathbb Z)$ is the standard torus.

{\smallskip\noindent}
We denote as ${\cal O}_{\rho, s, \delta, r, \xi}$ the set of complex holomorphic functions \begin{eqnarray*}\phi:\quad {\mathbb P}_{\hat\rho, \hat s, \hat\delta, \hat r, \hat\xi}\to {{\mathbb C}}\end{eqnarray*} for some  $\hat\rho>\rho$, $\hat s>s$, $\hat\delta>\delta$, $\hat r>r$, $\hat\xi>\xi$, equipped with 
 the norm
\begin{eqnarray*}\|\phi\|_{\rho, s, \delta, r, \xi}:=\sum_{k,h,j}\|\phi_{khj}\|_{\rho, r,\xi}e^{s|k|}\delta^{h+j}\end{eqnarray*}
where $\phi_{khj}({\mathbf I}, y,  x)$ are the coefficients of the Taylor--Fourier expansion\footnote{We denote as ${\mathbf x}^h:=x_1^{h_1}\cdots x_n^{h_n}$, where ${\mathbf x}=(x_1, \cdots, x_n)\in {{\mathbb R}}^n$ and $h=(h_1, \cdots, h_n)\in {{\mathbb N}}^n$.}
\begin{eqnarray*}\phi=\sum_{k,h,j}\phi_{khj}({\mathbf I},y, x)e^{{\rm i} k s}{\mathbf p}^h {\mathbf q}^j\ ,\quad \|\phi\|_{\rho, r,\xi}:=\sup_{{\mathbb I}_\rho\times{\mathbb Y}_r\times  {\mathbb X}_\xi}|\phi({\mathbf I}, y, x)|\ .\end{eqnarray*}
If $\phi$ is independent of $x$, we simply write
$\|\phi\|_{\rho, r}$ for $\|\phi\|_{\rho, r,\xi}$.

%{\smallskip\noindent}
%For a given vector--valued  function $\underline\phi=(\phi_1,\cdots, \phi_k)\in {\cal O}_{\rho, s, \delta, r, \xi}^k$, we let
%$$\|\underline\phi\|_{\rho, s, \delta, r, \xi}:=\sum_{i=1}^k \|\phi_i\|_{\rho, s, \delta, r, \xi}\ .$$

{\smallskip\noindent}
If $\phi\in {\cal O}_{\rho, s, \delta, r, \xi}$, we define its ``off--average'' $\widetilde\phi$ and   ``average'' $\overline\phi$ as
\begin{eqnarray*}&&\widetilde\phi:=\sum_{k,h,j:\atop (k,h-j)\ne (0,0)}\phi_{khj}({\mathbf I}, y,  x)e^{{\rm i} k s}{\mathbf p}^h {\mathbf q}^j\nonumber\\
&&\overline\phi:=\phi-\widetilde\phi=\frac{1}{2\pi^2}\int_{[0,2\pi]^n}\P_{{\mathbf p}\mathbf q}\phi({\mathbf I}, \bm\varphi, {\mathbf J}({\mathbf p}, {\mathbf q}), y, x)d\bm\varphi
\ ,\end{eqnarray*}
with
\begin{eqnarray*}\P_{{\mathbf p}\mathbf q}\phi({\mathbf I}, \bm\varphi, {\mathbf J}({\mathbf p}, {\mathbf q}), y, x):=\sum_{k,h}\phi_{khh}({\mathbf I}, y,  x)e^{{\rm i} k s}{\mathbf p}^h {\mathbf q}^h\end{eqnarray*}
We decompose
\begin{eqnarray*}{\cal O}_{\rho, s, \delta, r, \xi}={\cal Z}_{\rho, s, \delta, r, \xi}\oplus {\cal N}_{\rho, s, \delta, r, \xi}\ .\end{eqnarray*} where ${\cal Z}_{\rho, s, \delta, r, \xi}$, ${\cal N}_{\rho, s, \delta, r, \xi}$
are the ``zero--average'' and the  ``normal'' classes
\begin{eqnarray} \label{zero average}&&{\cal Z}_{\rho, s, \delta, r, \xi}:=\{\phi\in {\cal O}_{\rho, s, \delta, r, \xi}:\quad \phi=\widetilde\phi\}=\{\phi\in {\cal O}_{\rho, s, \delta, r, \xi}:\quad \overline\phi=0\}\\
\label{phi independent}&&{\cal N}_{\rho, s, \delta, r, \xi}:=\{\phi\in {\cal O}_{\rho, s, \delta, r, \xi}:\quad\phi=\overline\phi\}=\{\phi\in {\cal O}_{\rho, s, \delta, r, \xi}:\quad \widetilde\phi=0\}\ .\end{eqnarray} 
respectively. 
We finally let
$\omega_{y,{\mathbf I},{\mathbf J}}:=\partial_{y,{\mathbf I},{\mathbf J}} {\rm h}$.

{\smallskip\noindent}
We shall prove the following result, where no--resonance condition is required on the frequencies $\omega_{\mathbf I}$, which, as a matter of fact, might also be zero.

\begin{theorem}\label{NFL}
For any $n$, $m$, there exists a number ${\rm c}_{n,m}\ge 1$ such that, for any $N\in {\mathbb N}$ such that the following inequalities are satisfied
\begin{eqnarray} \label{normal form assumptions}
4N{\cal X}\left\|\Im\frac{\omega_{\mathbf I}}{\omega_y}\right\|_{\rho, r}<
s
\ ,\quad
4N{\cal X}\left\|\frac{\omega_{\mathbf J}}{\omega_y}\right\|_{\rho, r}<
1\ ,\quad   {\widetilde{\rm c}_{n,m}N\frac{{\cal X}}{{\rm d}}   % \|\frac{1}{\omega_y}\|_{\rho, r}
 \left\|{f}\right\|_{\rho, s, \delta, r, \xi}\left\|\frac{1}{\omega_y}\right\|_{\rho, s, \delta, r, \xi} <1}%\nonumber\\
  %&& {\rm c}^{\red2}_{n,m}N^{\red2}\frac{{\cal X}^{2}}{{\rm d}^{2}}  { \left\|\frac{1}{\omega_y}\right\|^2_{\rho, r} 
%\left\|% \frac
%{f}%{\omega_y}
%\right\|_{\rho, s, \delta, r, \xi}}\left\|% \frac
%\widetilde{f}%{\omega_y}
%\right\|_{\rho, s, \delta, r, \xi}<1 
\end{eqnarray} 
with ${\rm d}:=\min\big\{\rho s, r\xi, {\delta}^2\big\}$, ${\cal X}:=\sup\big\{|x|:\ x\in {\mathbb X}_\xi\big\}$, one can find an operator \begin{eqnarray} \label{PsiN}\Psi_*:\quad {\cal O}_{\rho, s, \delta, r, \xi}\to {\cal O}_{1/3 (\rho, s, \delta, r, \xi)}\end{eqnarray} 
which carries ${\rm H}$ to
\begin{eqnarray*}{\rm H}_*={\rm h}+g_*+f_*\end{eqnarray*}
where
$g_*\in {\cal N}_{1/3 (\rho, s, \delta, r, \xi)}$, $f_*\in {\cal O}_{1/3 (\rho, s, \delta, r, \xi)}$ and, moreover, the following inequalities hold 
\begin{eqnarray} \label{thesis}
&&\|g_*-\overline f\|_{1/3 (\rho, s, \delta, r, \xi)}\le 162\widetilde{\rm c}_{n, m} \frac{{\cal X}}{\rm d}\left\|%\frac{1}{\omega_y}\|_{\rho, r} \|
{\frac{\widetilde f}{\omega_y}}\right\|_{\rho, s, \delta, r, \xi}\| f\|_{\rho, s, \delta, r, \xi}\nonumber\\
&&    \|f_*\|_{1/3 (\rho, s, \delta, r, \xi)}\le \frac{1}{2^{N+1}} \|f\|_{\rho, s, \delta, r, \xi}\ .\end{eqnarray} 
The transformation $\Psi_*$ can be obtained as a composition of time--one Hamiltonian flows, and satisfies the following.  If
\begin{eqnarray*}({\mathbf I}, \bm\varphi, {\mathbf p}, {\mathbf q}, y, x):=\Psi_*({\mathbf I}_*, \bm\varphi_*, {\mathbf p}_*, {\mathbf q}_*, {\rm R}_*, {\rm r}_*)\end{eqnarray*} the following uniform bounds hold:
\begin{eqnarray} \label{phi close to id***}
&&{\rm d}\max\Big\{\frac{|{\mathbf I}-{\mathbf I}_*|}{\rho},\ \frac{|\bm\varphi -\bm\varphi _*|}{s},\ 
\frac{|{\mathbf p}-{\mathbf p}_*|}{\delta},\ \frac{ |{\mathbf q}-{\mathbf q}_*|}{\delta},\ \frac{|y-y_*|}{r},\frac{ |x-x_*|}{\xi} \Big\}\nonumber\\
&&\le\max\Big\{s|{\mathbf I}-{\mathbf I}_*|,\ \rho|\bm\varphi -\bm\varphi _*|,\ 
\delta |{\mathbf p}-{\mathbf p}_*|,\ \delta |{\mathbf q}-{\mathbf q}_*|,\ \xi |y-y_*|,\ r |x-x_*| \Big\}\nonumber\\
&&\leq%\yellow{4\, {\cal X}   % \|\frac{1}{\omega_y}\|_{\rho, r}
% \| \frac{f}{\omega_y}\|_{\rho, s, \delta, r, \xi}}
 {{19\, {\cal X}   % \|\frac{1}{\omega_y}\|_{\rho, r}
\left \| \frac{ f}{\omega_y}\right\|_{\rho, s, \delta, r, \xi}
}}%+
 %4\, \frac{{\cal X}^2 N}{{\rm d}\|\omega_y\|_{\rho, s, \delta, r, \xi}}  % \|\frac{1}{\omega_y}\|_{\rho, r}
 %\left\| \frac{\widetilde f}{\omega_y}\right\|_{\rho, s, \delta, r, \xi}
%{\|f\|_{\rho, s, \delta, r, \xi}}
 %}
 \ .
\end{eqnarray} 
\end{theorem}
\begin{remark}[Extensions]\rm
   \item[(i)] There is  an obvious  extension to the case that ${\mathbb I}_\rho$, ${\mathbb T}^n_s$ are replaced with $({\mathbb I}_1)_{\rho_1}\times \cdots\times ({\mathbb I}_n)_{\rho_n}$, ${\mathbb T}_{s_1}\times \cdots\times {\mathbb T}_{s_n}$. In this case, the number $s$ in the former equation in~(\ref{normal form assumptions})  is to be replaced with $\min_i\{s_i\}$. Moreover, the product $\rho\, s$ in the definition of  ${\rm d}$  is to be replaced with $\min_i\{\rho_i\,s_i\}$.  Finally, the bound in~(\ref{phi close to id***})  is to be changed taking into account the different sizes.
 \item[(ii)] If $f$ does not depend on some angle $\varphi_1$, $\cdots$, $\varphi_p$, the vector $\omega_{\mathbf I}$ in~(\ref{normal form assumptions}) is to be replaced with $\widehat\omega_{\mathbf I}:=(\omega_{{\mathbf I}_{p+1}},\cdots, \omega_{{\mathbf I}_{n}})$.
 %\item[(iii)]~(\ref{normal form assumptions}) $${\widetilde{\rm c}^2_{n,m}N^2\frac{{\cal X}^2}{{\rm d}^2}   % \|\frac{1}{\omega_y}\|_{\rho, r}
 %\left\|\frac{\widetilde f}{\omega_y}\right\|_{\rho, s, \delta, r, \xi}\left\|{ f}\right\|_{\rho, s, \delta, r, \xi}\left\|\frac{1}{\omega_y}\right\|_{\rho, s, \delta, r, \xi} <1}$$
\end{remark}

\subsection{Outline  of  the proof}
The complete proof of Theorem~\ref{NFL} is provided  in the Appendix~\ref{Time--one flows}, but here we aim spend some word,  so as to highlight the main  ideas.

{\smallskip\noindent}
We proceed by recursion.
We assume that, at a certain step, we have a system of the form
\begin{eqnarray} \label{step i}{\rm H}({\mathbf I}, \bm\varphi, {\mathbf J}({\mathbf p}, {\mathbf q}), y)={\rm h}({\mathbf I}, {\mathbf J}({\mathbf p}, {\mathbf q}), y)+g({\mathbf I}, {\mathbf J}({\mathbf p}, {\mathbf q}), y, x)+f({\mathbf I}, \bm\varphi, {\mathbf J}({\mathbf p}, {\mathbf q}), y, x)\end{eqnarray} 
where  $f\in {\cal O}_{\rho, s, \delta, r, \xi}$, $g\in {\cal N}_{\rho, s, \delta, r, \xi}$. At the first step, just take $g\equiv 0$.

{\smallskip\noindent}
After splitting $f$ on its Taylor--Fourier basis 
\begin{eqnarray*} f=\sum_{k,h,j} f_{khj}({\mathbf I},y, x)e^{{\rm i} k \cdot\bm\varphi}{\mathbf p}^h  {\mathbf q}^j\ .\end{eqnarray*}
one looks for a time--1 map 
\begin{eqnarray*}\Phi=e^{{\cal L}_\phi}=\sum_{k=0}^{\infty}\frac{{\cal L}_\phi^k}{k!}\qquad {\cal L}_\phi(f):=\big\{\phi,\ f\big\}\end{eqnarray*} generated by a small Hamiltonian 
$\phi$
which will be taken  in the  class ${\cal Z}_{\rho, s, \delta, r, \xi}$ in~(\ref{zero average}).
Here,
\begin{eqnarray*} \big\{\phi,\ f\big\}:=\sum_{i=1}^n(\partial_{{\mathbf I}_i}\phi\partial_{\bm\varphi_i}f-\partial_{{\mathbf I}_i}f\partial_{\bm\varphi_i}\phi)+\sum_{i=1}^m(\partial_{{\mathbf p}_i}\phi\partial_{\mathbf q_i}f-\partial_{{\mathbf p}_i}f\partial_{\mathbf q_i}\phi)+(\partial_{y}\phi\partial_{x}f-\partial_{y}f\partial_{x}\phi)\end{eqnarray*}
denotes the Poisson parentheses of $\phi$ and $f$.
One lets
\begin{eqnarray} \label{exp}\phi=\sum_{(k,h,j):\atop{(k,h-j)\ne (0,0)}} \phi_{khj}({\mathbf I},y, x)e^{{\rm i} k\cdot\bm\varphi}{\mathbf p}^h {\mathbf q}^j\ .\end{eqnarray} 

{\smallskip\noindent}
The operation
\begin{eqnarray*} \phi\to \{\phi,{\rm h}\}\end{eqnarray*}
acts diagonally  on the monomials  in the expansion~(\ref{exp}), carrying
\begin{eqnarray} \label{diagonal}\phi_{khj}\to -\big(\omega_y\partial_x \phi_{khj}+\l_{khj} \phi_{khj}\big)\ ,\quad {\rm with}\quad \l_{khj}:=(h-j)\cdot\omega_{\mathbf J}+{\rm i} k\cdot\omega_{\mathbf I}\ .\end{eqnarray} 
Therefore, one defines
\begin{eqnarray*}\{\phi,{\rm h}\}=:-D_\omega\phi\ .\end{eqnarray*}
The formal application of $\Phi=e^{{\cal L}_\phi}$ yields:
\begin{eqnarray} \label{f1}
e^{{\cal L}_\phi} {\rm H}&=&e^{{\cal L}_\phi} ({\rm h}+g+f)={\rm h}+g-D_\omega \phi+f+\Phi_2({\rm h})+\Phi_1(g)+\Phi_1(f)
\end{eqnarray} 
where the $\Phi_h$'s are the tails of $e^{{\cal L}_\phi}$, defined in Appendix~\ref{Time--one flows}.

{\smallskip\noindent}
Next, one requires that   the residual term $-D_\omega \phi+f$ lies in the class ${\cal N}_{\rho, s, \delta, r, \xi}$ in~(\ref{phi independent})%. This
%amounts  to solve 
%the ``homological'' equation
\begin{eqnarray} \label{homological equation}{(-D_\omega \phi+f
)}\in {\cal N}_{\rho, s, \delta, r, \xi}\end{eqnarray} 
for $\phi$.

{\smallskip\noindent}
Since we have chosen $\phi\in{\cal Z}_{\rho, s, \delta, r, \xi}$, by~(\ref{diagonal}), we have that also $D_\omega \phi\in{\cal Z}_{\rho, s, \delta, r, \xi}$. So, Equation~(\ref{homological equation}) becomes
\begin{eqnarray*}-D_\omega \phi+\widetilde f=0\ .\end{eqnarray*}
In terms of  the Taylor--Fourier modes, the equation becomes
\begin{eqnarray} \label{homol eq} \omega_y\partial_x \phi_{khj}+\l_{khj} \phi_{khj}=f_{khj}\qquad \forall\ (k,h,j):\ (k,h-j)\ne (0,0)\ .\end{eqnarray}

{\smallskip\noindent} In the standard situation, one typically proceeds to solve such equation via Fourier series:
\begin{eqnarray*} f_{khj}({\mathbf I},y, x)=\sum_{\ell}f_{khj\ell}({\mathbf I}, y)e^{{\rm i} \ell x}\ ,\qquad \phi_{khj}({\mathbf I},y, x)=\sum_{\ell}\phi_{khj\ell}({\mathbf I}, y)e^{{\rm i} \ell x}\end{eqnarray*}
so as to find
$\displaystyle\phi_{khj\ell}=\frac{f_{khj\ell}}{\mu_{khj\ell}}$
with the usual  denominators $\mu_{khj\ell}:=\l_{khj}+{\rm i}\ell \omega_y$ which one requires not to vanish %(equivalently, $\frac{\l_{khj}}{\omega_y}\notin {\rm i}{\mathbb Z}$)
via, e.g., a ``diophantine inequality'' to be held for all $(k,h,j,\ell)$ with $(k,h-j)\ne (0,0)$. In this standard case, there is not much freedom in the choice of $\phi$. In fact, such solution is determined up to solutions  of the homogenous equation
\begin{eqnarray} \label{homogeneous}D_\omega\phi_0=0\end{eqnarray} 
which, in view of the Diophantine condition, has the only trivial solution $\phi_0\equiv0$. {\it The situation is different if $f$ is not periodic in $x$, or $\phi$ is not needed so}. In such a case, it is possible to find a solution of~(\ref{homol eq}), corresponding to a non--trivial solution of~(\ref{homogeneous}),  where small divisors do not appear. This is
\begin{eqnarray} \label{solution}\phi_{khj}({\mathbf I}, y,  x)=\left\{\begin{array}{llll}\displaystyle
\frac{1}{\omega_y}\int_0^xf_{khj}({\mathbf I}, y, \tau)e^{\frac{\l_{khj}}{\omega_y}(\tau-x)}d\tau\quad {\rm if}\ \ (k,h-j)\ne (0,0)\\\\
\displaystyle0\qquad\qquad\qquad\qquad\qquad\qquad\qquad\ {\rm otherwise .}
\end{array}\right.
\end{eqnarray} 
Multiplying by $e^{ik\varphi}$ and summing over $k$, $h$ and $j$, we obtain
\begin{eqnarray} \label{phiNEW}\phi({\mathbf I}, {\bm\varphi}, p, q, y, x)=\frac{1}{\omega_y}\int_0^x \widetilde f\left({\mathbf I}, {\bm\varphi}+\frac{\omega_{\mathbf I}}{\omega_y}(\tau-x), p e^{\frac{\omega_{\mathbf J}}{\omega_y}(\tau-x)}, q e^{-\frac{\omega_{\mathbf J}}{\omega_y}(\tau-x)}, y, \tau\right)d\tau\ .\end{eqnarray} 
In Appendix~\ref{Time--one flows}, we shall prove that, under the assumptions~(\ref{normal form assumptions}), this 
function can be used
to obtain a convergent time--one map and that the construction can be iterated so as to provide the proof of Theorem~\ref{NFL}. The construction of the iterations and the proof of its convergence is obtained adapting the techniques of
~\cite{poschel93} to the present case.

\section{Proof of Theorem~\ref{main}}\label{Completion of the proof of Theorem}
%with
%$$\widehat{\rm V}_\varepsilon({\rm G}, {\rm g}):=\widehat{\rm U}_\varepsilon({\rm G}, {\rm g})-1\ .$$
The   purpose of this section is state and prove a more precise version  of Theorem~\ref{main} (Theorem~\ref{main1} below). We shall obtain it as an 
application of Theorem~\ref{NFL} to the Hamiltonians $\widehat{\rm H}_i$. Therefore,    we need to introduce a change new coordinates ${\cal C}$
which put the $\widehat{\rm H}_i$'s in the suited form~(\ref{h0f0}). 
Since the potentials ${\rm V}_i$ in~(\ref{Vi}) are, for large ${\rm r}$, asymptotic to $-\frac{m_0^2}{{\rm r}}$, it is convenient to rewrite the functions $\widehat{\rm H}_i$ in~(\ref{VVV}) as
\begin{eqnarray*}
\widehat{\rm H}_1({\rm R}, {\rm G}, {\rm r}, {\rm g})&=&\left(\frac{{\rm R}^2}{2m_0}-\frac{m_0^2}{{\rm r}}\right)+\frac{{\rm G}^2}{2m_0{\rm r}^2}-\frac{\overline\beta}{\beta+\overline\beta}\frac{m_0^2}{{\rm r}}\left(\widehat{\rm F}_{-\overline\beta\varepsilon({\rm r})}\Big(\widehat{\rm E}_{\beta\varepsilon({\rm r})}({\rm G}, {\rm g})\Big)-1\right)\nonumber\\\nonumber\\
&&-\frac{\beta}{\beta+\overline\beta}\frac{m_0^2}{{\rm r}}\left(\widehat{\rm F}_{-\overline\beta\varepsilon({\rm r})}\Big(\widehat{\rm E}_{-\overline\beta\varepsilon({\rm r})}({\rm G}, {\rm g})\Big)-1\right)\nonumber\\
\widehat{\rm H}_2({\rm R}, {\rm G}, {\rm r}, {\rm g})&=&\left(\frac{{\rm R}^2}{2m_0}-\frac{m_0^2}{{\rm r}}\right)+\frac{{\rm G}^2}{2m_0{\rm r}^2}-\frac{\overline\beta}{\beta+\overline\beta}\frac{m_0^2}{{\rm r}}\left(\widehat{\rm F}_{(\beta+\overline\beta)\varepsilon({\rm r})}\Big(\widehat{\rm E}_{-\overline\beta\varepsilon({\rm r})}({\rm G}, {\rm g})\Big)-1\right)
\end{eqnarray*}
and take ${\cal C}$ as the composition of two independent and canonical changes
\begin{eqnarray*}{\cal C}_1:\ ({\cal G}, \gamma)\to ({\rm G}, {\rm g})\ ,\quad {\cal C}_2:\ (y, x)\to ({\rm R}, {\rm r})\end{eqnarray*} where 
${\cal C}_1$ is defined as in~(\ref{Gg}), with $k=0$, while ${\cal C}_2$ is defined via the formulae
 \begin{eqnarray} \label{Rr}\left\{\begin{array}{llll}\displaystyle{\rm R}(y, x)=\frac{m_0^3}{y}\sqrt{ \frac{\cos\xi'(x)+1}{1-\cos\xi'(x)}}
 \\\\
 \displaystyle {\rm r}(y, x)=\frac{y^2}{m_0^3}(1-\cos\xi'(x))\end{array}\right.\end{eqnarray} 
where $\xi'(x)$ solves
\begin{eqnarray} \label{Kepler1}\xi'-\sin\xi'=  x\ .\end{eqnarray} 
${\cal C}_2$ has been chosen so that
\begin{eqnarray*} {\rm H}_\omega\circ{\cal C}_1=\left(\frac{{\rm R}^2}{2m_0}-\frac{m_0^2}{{\rm r}}\right)\circ{\cal C}_1=-\frac{m_0^5}{2y^2}\ .\end{eqnarray*} Using the new coordinates, we have
\begin{eqnarray*}
\widehat{\rm H}_1&=&-\frac{m_0^5}{2y^2}+\frac{m_0^2}{{\rm r}(y, x)}\left(\varepsilon(y, x)\frac{({\Lambda}^2-{\cal G}^2)}{2{\Lambda}^2}\cos^2\gamma-\frac{\overline\beta}{\beta+\overline\beta}\left(\widehat{\rm F}%^+
_{\beta\varepsilon(y, x)}\Big(\widehat{\rm E}_{\beta\varepsilon(y, x)}({\cal G}, \gamma)\Big)-1\right)\right.\nonumber\\
&&\left.-\frac{\beta}{\beta+\overline\beta}\left(\widehat{\rm F}_{-\overline\beta\varepsilon(y, x)}\Big(\widehat{\rm E}_{-\overline\beta\varepsilon(y, x)}({\cal G}, \gamma)\Big)-1\right)\right)\nonumber\\\nonumber\\
\widehat{\rm H}_2&=&-\frac{m_0^5}{2y^2}+\frac{m_0^2}{{\rm r}(y, x)}\left(\varepsilon(y, x)\frac{({\Lambda}^2-{\cal G}^2)}{2{\Lambda}^2}\cos^2\gamma\right.\nonumber\\
&&\left.
-
\frac{\overline\beta}{\beta+\overline\beta}\frac{m_0^2}{{\rm r}(y, x)}\left(\widehat{\rm F}%^+
_{(\beta+\overline\beta)\varepsilon(y, x)}\Big(\widehat{\rm E}_{(\beta+\overline\beta)\varepsilon(y, x)}({\cal G}, \gamma)\Big)-1\right)\right)
\end{eqnarray*}
having abusively denoted as $\varepsilon(y, x)$ the function $\varepsilon({\rm r}(y, x))$, and (using the formulae in~(\ref{EActionAngle})--(\ref{UUU})
\begin{eqnarray} \label{hatE}\widehat{\rm E}_\varepsilon({\cal G}, \gamma):=\frac{{\cal G}}{\Lambda}+\varepsilon\left(1-\frac{{\cal G}^2}{{\Lambda}^2}\right)\cos^2\gamma\ .\end{eqnarray} 
We now determine a domain of holomorphy of $\widehat{\rm H}_i$. 
{Recall that we use, as mass parameters, the numbers $\beta_*$, $\beta^*$ in~(\ref{b*})}. For the coordinates $(y, x)$, we choose the {complex} domains ${\mathbb Y}_{\sqrt{m_0^3\alpha_-}}$, ${\mathbb X}_{\sqrt{\varepsilon_0}}$, where
\begin{eqnarray} \label{yx}{\mathbb Y}:=\Big\{y\in {\mathbb R}:\ 2\sqrt{m_0^3\alpha_-}<y<\sqrt{m_0^3\alpha_+}\Big\}\ ,\quad {\mathbb X}:=\Big\{x\in {\mathbb R}:\ |x-\pi|\le \pi-2\sqrt{\varepsilon_0}\Big\}\end{eqnarray} 
with $0<\varepsilon_0<1$, $0<\alpha_-<\alpha_+{/4}$ verifying
\begin{eqnarray} \label{Kepineq}\alpha_-\varepsilon_0>\frac{4\beta^* a}{c_0}\qquad {\rm with}\end{eqnarray} 
with $\beta^*$ as in~(\ref{b*}) and $c_0$ being such that   for any $0<\varepsilon_0<1$ and for any  $x\in{\cal X}_{\sqrt{\varepsilon_0}}$,
Equation~(\ref{Kepler1}) has a unique solution $\xi'(x)$ which depends analytically on $x$ and and verifies \begin{eqnarray} \label{xiineq} |1-\cos \xi'(x)|\ge c_0\varepsilon_0\ .\end{eqnarray} 
(the existence of such a number $c_0$  is well known).
For the coordinates ${\cal G}$, $\gamma$, we choose the domains
$ {\mathbb G}_\delta$, ${\mathbb T}_{s_0}$,
with $0<\delta<\Lambda$, $s_0$ fixed, and ${\mathbb G}:=\Big\{{\cal G}\in {\mathbb R}:\ \Lambda-\delta<{\cal G}<\Lambda\Big\}$. Remark that, since, in $\widehat{\rm H}_i$, there is no dependence of ${\rm G}$, but only on ${\rm G}^2$, ${\cal G}=\Lambda$ is a regular point for $\widehat{\rm H}_i$.
Then we let
\begin{eqnarray*} {\mathbb D}:= {\mathbb G}\times {\mathbb T}\times {\mathbb Y}\times {\mathbb X}\subset {\mathbb R}^4\end{eqnarray*}
and
\begin{eqnarray*}{\mathbb D}_{\delta, s_0, \sqrt{m_0^3\alpha_-},\sqrt{\varepsilon_0}}:= {\mathbb G}_\delta\times {\mathbb T}_{s_0}\times {\mathbb Y}_{\sqrt{m_0^3\alpha_-}}\times {\mathbb X}_{\sqrt{\varepsilon_0}}\subset {\mathbb C}^4\ .\end{eqnarray*}

{\smallskip\noindent}
{We now check that, under the further assumptions
\begin{eqnarray} \label{ass1}0<\delta\le \frac{\Lambda}{4}\ ,\qquad C^*(s_0)\frac{\delta}{\Lambda}<1\qquad C^*(s_0):=16\left(\sup_{ {\mathbb T}_{s_0}}|\sin\gamma|\right)^2\end{eqnarray} 
}
   $\widehat{\rm H}_i$ are holomorphic functions on the domain 
\begin{eqnarray} \label{domain}{\mathbb D}_{\delta, s_0, \sqrt{m_0^3\alpha_-},\sqrt{\varepsilon_0}}:={\mathbb Y}_{\sqrt{m_0^3\alpha_-}}\times {\mathbb X}_{\sqrt{\varepsilon_0}}\times {\mathbb G}_\delta\times {\mathbb T}_{s_0}\ .\end{eqnarray} 
By 
(\ref{xiineq}),  the first equation in~(\ref{yx}), and the expression of ${\rm r}(y, x)$ in~(\ref{Rr}), we have
\begin{eqnarray} \label{rineq}|{\rm r}(y, x)|\ge c_0\alpha_-\varepsilon_0 \end{eqnarray} 
and hence, because of~(\ref{Kepineq}),
\begin{eqnarray} \label{epsineq}|\beta^*\varepsilon(y, x)|=\left|\frac{\beta^*a}{{\rm r}(y, x)}\right|\le\frac{\beta^* a}{c_0\alpha_-\varepsilon_0}<\frac{1}{4}\ . \end{eqnarray} 
{By inequalities~(\ref{rineq})--(\ref{epsineq}) and Proposition~\ref{holomorphy}, we only need to check that, if  $\varepsilon_*\in\{\beta\varepsilon, -\overline\beta\varepsilon\}$ for $i=1$ and $\varepsilon_*=(\beta+\overline\beta)\varepsilon$ for $i=2$, then  Equation~(\ref{sing}) with $\varepsilon=\varepsilon_*$ and $t=\widehat{\rm E}_{\varepsilon_*}({\cal G}, \gamma)$ has not solutions in ${\mathbb D}_{\delta, s_0, \sqrt{m_0^3\alpha_-},\sqrt{\varepsilon_0}}$.  
We prove that any such solution would verify $|\varepsilon_*|\ge \frac{1}4$,  which would contradict~(\ref{epsineq}), as
$|\beta^*\varepsilon|\ge |\varepsilon_*|$. 
 Using~(\ref{hatE}), Equation~(\ref{sing}) with $\varepsilon=\varepsilon_*$ and $t=\widehat{\rm E}_{\varepsilon_*}({\cal G}, \gamma)$ is
\begin{eqnarray*} 4\varepsilon_*^2\left(1-\left(1-\frac{{\cal G}^2}{{\Lambda}^2}\right)\cos^2\gamma\right)-4\frac{{\cal G}}{\Lambda}\varepsilon_*+1=0\ .\end{eqnarray*} We solve for $\varepsilon_*$:
\begin{eqnarray*} \varepsilon_*=\frac{1}2\left(\frac{{\cal G}}{\Lambda}+\sin\gamma\sqrt{\frac{{\cal G}^2}{{\Lambda}^2}-1}\right)\end{eqnarray*}
with the double determination of the square root. Since $\left|\frac{{\cal G}}{\Lambda}\right|\ge 1-\frac{\delta}{\Lambda}\ge \frac{3}{4}$ and, if $c_*(s_0):=\sup_{{\mathbb T}_{s_0}}|\sin\gamma|$, $\left|\sin\gamma\sqrt{\frac{{\cal G}^2}{{\Lambda}^2}-1}\right|\le c^*(s_0)\sqrt{\frac{\delta}{\Lambda}}\le \frac{1}4$,  we have $|\varepsilon_*|\ge \frac{1}2\left(\frac{3}{4}-\frac{1}4\right)=\frac{1}4$, as claimed.
 }

{\smallskip\noindent}
We are now ready to state the result.
Observe  that the motions 
\begin{eqnarray*} {\cal G}={\cal G}_0={\rm constant}\ ,\quad |\gamma(T)-\gamma(0)|=2\pi\end{eqnarray*}
correspond, using the coordinates $({\rm G}, {\rm g})$, to librations about $(0,0)$ if $0<{\cal G}_0<\Lambda$; about $(0, \pi)$ if $-\Lambda<{\cal G}_0<0$.  We shall prove the existence, in $\widehat{\rm H}_i$'s, of motions close to these ones.

\begin{theorem}\label{main1}
Let 
$\alpha_-$, $\alpha_+$, $\beta$, $\overline\beta$, $\delta$, $\varepsilon_0$ and $s_0$ %, ${\mathbb D}_{\d, s_0, \sqrt{m_0^3\a_-},\sqrt{\varepsilon_0}}$
 be fixed; $\beta_*$, $\beta^*$ as in ~(\ref{b*}).
There exist two numbers  $C^*>C_*>1$, both independent of $\alpha_-$, $\alpha_+$, $\beta$, $\overline\beta$, $\delta$, $\varepsilon_0$, with $C^*$ possibly depending on  $s_0$, while  $C_*$ independent of  $s_0$, such that, if the following inequalities are satisfied
 \begin{eqnarray} \label{N0}
 &&0<\varepsilon_0<1\ ,\quad 0<\delta\le\frac{\Lambda}{4}\ ,\quad \frac{4\beta^* a}{c_0\alpha_-\varepsilon_0}<1\ ,\quad \frac{C^*\delta}{\beta_*\Lambda}\le 1\nonumber\\
 &&\frac{1}{N_0}:=C_*\max
\left\{%\frac{1}{c_0^2\varepsilon_0^2}\sqrt{\frac{a}{\alpha_-}}\ ,\quad
%\frac{\delta}{c_0^2\varepsilon_0^2\sqrt{m_0^3\alpha_-\varepsilon_0}}\sqrt{\frac{a}{\alpha_-}}\ ,\quad
\frac{\beta_*\Lambda}{c_0^2\varepsilon_0^2\delta s_0}\sqrt{\frac{a}{\alpha_-}}\	 ,%\right.\nonumber\\
%&&\left.	
	\quad
 {\frac{\beta_*}{c_0^2\varepsilon_0^{5/2}}{\frac{a}{\alpha_-}}}
\right\}\frac{\alpha_+^{3/2}}{\alpha_-^{3/2}}<\frac{c_0^2\varepsilon_0^2\alpha_-^2}{2\alpha_+^2}
\end{eqnarray} 
it is possible to find new coordinates $({\cal G}_*, \gamma_*, y_*, x_*)$ and a time $T$ such that any solution $\Gamma_*(t)=({\cal G}_*(t), \gamma_*(t), y_*(t), x_*(t)) $ of $\widehat{\rm H}_i$ with initial datum $\Gamma_*(0)$ $=$ $({\cal G}_*(0)$, $\gamma_*(0)$, $y_*(0)$, $x_*(0))$ $\in{\mathbb D}$ such that 
\begin{eqnarray} \label{initial data}|{\cal G}(0)-\Lambda|\le \frac{\delta}{2}\ ,\quad 2\sqrt{m_0^3\alpha_-}\le|y_*(0)|{\le\frac{\sqrt{m_0^3\alpha_-}+\sqrt{m_0^3\alpha_+}}{2} } \ ,\quad x_*(0)=\pi
\end{eqnarray} 
stays in ${\mathbb D}$ for all $0\le t\le T$ and ${\cal G}_*$ varies a little in the course of such time:
\begin{eqnarray*}|{\cal G}_*(t)-{\cal G}_*(0)|\le C_*\frac{2^{-N_0}}{s_0}\frac{m_0^2a\beta_*}{c^2_0\varepsilon^2_0\alpha^2_-}t \ \quad {\rm for\ all}\ 0\le t\le T\ .\end{eqnarray*}
If, in addition,
\begin{eqnarray} \label{last cond}{\eta}:=C_*\max\left\{\frac{\alpha_+^2}{\beta_*\sqrt{\alpha_-^3 a}}\ ,\ \frac{\alpha_+^2}{c^2_0\varepsilon^{5/2}_0\alpha^2_-}\sqrt{\frac{a}{\alpha_-}}\ ,\ 
\frac{\alpha_+^2}{c^2_0\varepsilon^{2}_0\alpha^2_-} \frac{\Lambda}{s_0\delta}2^{-N_0}
\right\}<1
\end{eqnarray} 
then motions close to librations occur, in the sense that, also
\begin{eqnarray*}|\gamma_*(T)-\gamma_*(0)|\ge 3\pi\ .\end{eqnarray*}
The time $T$ can be taken to be
	\begin{eqnarray} \label{T}
	{T=\frac{\Lambda\alpha_+^3}{\beta_*m_0^2a}\frac{3\pi}{\eta}}
	\end{eqnarray} 
Finally, the change
\begin{eqnarray*}({\cal G}_*, \gamma_*, y_*, x_*)\to ({\cal G}, \gamma, y, x)\end{eqnarray*}
is real--analytic and close to the identity, in the sense that
\begin{eqnarray*}|{\cal G}-{\cal G}_*|\le \frac{\Lambda}{N_0}\ ,\quad |\gamma-\gamma_*|\le \frac{s_0}{N_0}\ ,\quad |y-y_*|\le \frac{\sqrt{m_0^3\alpha_-}}{N_0}\ ,\quad |x-x_*|\le \frac{\sqrt{\varepsilon_0}}{N_0}\ .\end{eqnarray*}
\end{theorem}
{\begin{remark}[Proof of Theorem~\ref{main}]\label{main1rk}\rm Inequalities~(\ref{N0}),~(\ref{initial data}) and~(\ref{last cond}) are simultaneously satisfied provided that the following holds. Fix $0<\varepsilon_0<1$ and $0<\delta\le \frac{\Lambda}{4}$ once forever. Then identify $\delta$ as the size of ${\rm U}_0$ and $2^{-N_0}$ as the size of ${\rm V}_0$. Take
\begin{eqnarray*}
&&s_0\rhd\frac{\Lambda}{\delta\varepsilon_0^4}\ ,\quad \frac{\alpha_+}{a}=2^8\frac{\alpha_-}{a}\rhd\frac{{C^*(s_0)}^2}{\varepsilon_0^8}\nonumber\\
&&\min\left\{
\varepsilon_0^4\frac{\delta}{\Lambda}s_0\sqrt{\frac{\alpha_-}{a}},\ \varepsilon_0^{9/2}\frac{\alpha_-}{a}
\right\}\rhd\beta^*\ge \beta_*\rhd\max\left\{C^*\frac{\delta}{\Lambda},\ \sqrt{\frac{\alpha_-}{a}}\right\}
\end{eqnarray*} 
For short, we have written ``$a\rhd b$'' if there exist $c>1$, independent of $\delta$, $\Lambda$, $\alpha_-$, $\alpha_+$, $a$ and $s_0$ such that $a>cb$.
Note that here it is essential that $\alpha_-$, $s_0$, $\beta_*$ and $\beta^*$ can be chosen arbitrarily large.
\end{remark}}

\par\medskip\noindent{\bf Proof\ } 
During the proof we shall make extensive use of Cauchy\footnote{{Observe that the assumptions~(\ref{N0}) include, in particular,~(\ref{Kepineq}) and~(\ref{ass1}), so $\widehat{\rm H}_i$'s are holomorphic on the domain~(\ref{domain}). }} inequalities. \\
We  aim to apply Theorem~\ref{NFL}, with ${\mathbf I}={\cal G}$, $\bm\varphi=\gamma$,
$(y, x)$ as in~(\ref{Rr}),
${\rm h}(y)=-\frac{m_0^5}{2y^2}$
and, finally

\begin{eqnarray} \label{f0}
f({\cal G}, \gamma, y, x)=\left\{\begin{array}{llll}\frac{m_0^2}{{\rm r}(y, x)}\left(\varepsilon(y, x)\frac{({\Lambda}^2-{\cal G}^2)}{2{\Lambda}^2}\cos^2\gamma-\frac{\overline\beta}{\beta+\overline\beta}\left(\widehat{\rm F}%^+
_{\beta\varepsilon(y, x)}\Big(\widehat{\rm E}_{\beta\varepsilon(y, x)}({\cal G}, \gamma)\Big)-1\right)\right.\\
\left.-\frac{\beta}{\beta+\overline\beta}\left(\widehat{\rm F}_{-\overline\beta\varepsilon(y, x)}\Big(\widehat{\rm E}_{-\overline\beta\varepsilon(y, x)}({\cal G}, \gamma)\Big)-1\right)\right)\quad i=1\\\\
\frac{m_0^2}{{\rm r}(y, x)}\left(\varepsilon(y, x)\frac{({\Lambda}^2-{\cal G}^2)}{2{\Lambda}^2}\cos^2\gamma\right.\\
\left.-\frac{\overline\beta}{\beta+\overline\beta}\frac{m_0^2}{{\rm r}(y, x)}\left(\widehat{\rm F}%^+
_{(\beta+\overline\beta)\varepsilon(y, x)}\Big(\widehat{\rm E}_{(\beta+\overline\beta)\varepsilon(y, x)}({\cal G}, \gamma)\Big)-1\right)\right)\\
\quad i=2
\end{array}\right.
\end{eqnarray} 
In our case, ${\mathbf p}$, ${\mathbf q}$ do not exist and the unperturbed term ${\rm h}$ does not depend on ${\mathbf I}={\cal G}$. Therefore, we have only to verify the last condition in~(\ref{normal form assumptions}). We have
\begin{eqnarray*} \omega_y=\frac{m_0^5}{y^3}\ ,\quad {\rm d}=\min\{\delta s_0\ ,\ \sqrt{m_0^3\alpha_-\varepsilon_0}\}\ ,\quad {\cal X}=\sqrt{4\pi^2+\varepsilon_0}\le 3\pi\end{eqnarray*}
(having used $\varepsilon_0<1$) and  
\begin{eqnarray} \label{Delta2}
&&{\left\|\frac{1}{\omega_y}\right\|_{\sqrt{m_0^3\alpha_-}, \sqrt{\varepsilon_0}, \delta, s_0}\le2\sqrt{\frac{\alpha_+^3}{m_0}}}\nonumber\\
&&{\left.\|f\|\right|_{\sqrt{m_0^3\alpha_-}, \sqrt{\varepsilon_0}, \delta, s_0}\le\frac{m_0^2a}{c^2_0\varepsilon^2_0\alpha^2_-}\left(
C_1\frac{\delta}{\Lambda}+C_2\beta_*
\right)\le C_*\frac{m_0^2a\beta_*}{c^2_0\varepsilon^2_0\alpha^2_-}=:\Delta}
\end{eqnarray} 
 with\footnote{
Split $f$ in~(\ref{f0}) as
$f=f_1+f_2$, where $f_1=\frac{m_0^2}{{\rm r}(y, x)}\varepsilon(y, x)\frac{({\Lambda}^2-{\cal G}^2)}{2{\Lambda}^2}\cos^2\gamma$.
The term proportional to $C_i$ corresponds to be a upper bound of  $\left\|\frac{f_i}{\omega_y}\right\|$.  $C_2$ can be chosen to be independent of $s_0$ because $\|f_2\|$ goes to zero as $s_0\to \infty$.
}  $C_1$, $C_2$, {$C_*$} independent of $\alpha_-$, $\alpha_+$, $\delta$,  {$\beta_*$, $\beta^*$} but $C_1$ possibly depending on $s_0$  while  $C_2$, {$C_*$} independent  of $s_0$.
We {have} choosen the number $C^*$ in~(\ref{N0}) larger or equal than $2C_1/C_2$ {and the number $C_*$ larger or equal than $3C_2/2$}, so that
$\left(C_1\frac{\delta}{\Lambda} +C_2\beta_*\right)\le \frac{3}{2}C_2\beta_*{\le C_*\beta_*}$.
 we have
\begin{eqnarray*}
\widetilde{\rm c}_{1, 0}\frac{\chi}{{\rm d}}\big\|f\big\|_{\sqrt{m_0^3\alpha_-}, \sqrt{\varepsilon_0}, \delta, s_0}\left\|\frac{1}{\omega_y}\right\|_{\sqrt{m_0^3\alpha_-}, \sqrt{\varepsilon_0}, \delta, s_0}&\le&  C_*\max
\left\{
\frac{\beta_*\Lambda}{c_0^2\varepsilon_0^2\delta s_0}\sqrt{\frac{a}{\alpha_-}}\	 ,		\right.\nonumber\\
&&\left.
\frac{\beta_*\Lambda}{c_0^2\varepsilon_0^2\sqrt{m_0^3\alpha_-\varepsilon_0}}\sqrt{\frac{a}{\alpha_-}}
\right\}\frac{\alpha_+^{3/2}}{\alpha_-^{3/2}}\nonumber\\
&&{
=C_*\max
\left\{
\frac{\beta_*\Lambda}{c_0^2\varepsilon_0^2\delta s_0}\sqrt{\frac{a}{\alpha_-}}\	 ,		\right.}\nonumber\\
&&{\left.
\frac{\beta_*}{c_0^2\varepsilon_0^{\frac{5}{2}}}{\frac{a}{\alpha_-}}
\right\}\frac{\alpha_+^{3/2}}{\alpha_-^{3/2}}
}\nonumber\\
&\le&\frac{1}{N_0}
\end{eqnarray*}
Therefore, the last condition in~(\ref{normal form assumptions}) is immediately implied by $N<N_0$, with $N_0$ as in~(\ref{N0}). We then find a real--analytic transformation 
\begin{eqnarray*}
\phi_*:\quad ({\cal G}_*, \gamma_*, y_*, x_*)\in{\mathbb D}_{\sqrt{m_0^3\alpha_-}/3,\sqrt{\varepsilon_0}/3, \delta/3, s_0/3}\to ({\cal G}, \gamma, y, x)\in{\mathbb D}_{\delta, s_0, \sqrt{m_0^3\alpha_-},\sqrt{\varepsilon_0}}\end{eqnarray*}
which leads $\widehat{\rm H}_i$ to
\begin{eqnarray} \label{stability}\widehat{\rm H}_{*}={\rm h}(y_*)+g_{*}(y_*, x_*, {\cal G}_*)+f_{*}({\cal G}_*, \gamma_*, y_*, x_*)\end{eqnarray} 
where, $g_*$, $f_*$ satisfy the  following bounds:
\begin{eqnarray*}\|g_*{-\overline f}\|\le 2 \Delta\ ,\quad \|g_*\|\le 2^{-N}\Delta \end{eqnarray*}
with {$\overline f(y_*, x_*, {\cal G}_*)$ the $\gamma_*$--average of $f(y_*, x_*, {\cal G}_*, \gamma_*)$ and $\Delta$ as in~(\ref{Delta2}).}
Let now $\Gamma_*(t)=({\cal G}_*(t), \gamma_*(t), y_*(t), x_*(t)) $ be a solution of $\widehat{\rm H}_i$ with initial datum $\Gamma_*(0)$ $=$ $({\cal G}_*(0)$, $\gamma_*(0)$, $y_*(0)$, $x_*(0))$ $ \in{\mathbb D}$ and verifying~(\ref{initial data}).
We look for a time $T>0$ such that $\Gamma_*(t)\in{\mathbb D}$
for all  $0\le t\le T$.  
We show that we can take $T$ as in~(\ref{T}), {which, for convenience we rewrite as
	\begin{eqnarray} \label{TOLD}T=\min\left\{\sqrt{\frac{\alpha_-^3}{m_0}}\ ,\ {\frac{\sqrt{m_0^3\alpha_-\varepsilon_0}}{\Delta}}\ ,\ 2^{N_0}\frac{s_0\delta}{\Delta}\right\}\end{eqnarray} 
	where $\Delta$ is as in~(\ref{Delta2}).}
	Equation~(\ref{stability}) implies
\begin{eqnarray*}
\left|y_*(t)-y_*(0)\right|\le \frac{\Delta t}{\sqrt{\varepsilon_0}}\end{eqnarray*}
So, for $t\le \frac{\sqrt{m_0^3\alpha_-\varepsilon_0}}{\Delta}$, we have
\begin{eqnarray*} t  \le \frac{|y_*(0)|-\sqrt{m_0^3\alpha_-}}{\Delta}\sqrt{\varepsilon_0}\quad \Longrightarrow\quad |y_*(t)-y_*(0)|\le |y_*(0)|-\sqrt{m_0^3\alpha_-} \end{eqnarray*}
namely, $y_*(t)\in {\mathbb Y}$ for $0\le t\le T$.
Since $|y|\ge \sqrt{m_0^3 \alpha_-}$ for all this time, we also have
\begin{eqnarray*}\left|x_*(t)-x_*(0)\right|\le \left(\sqrt{\frac{m_0}{\alpha_-^3}}+\frac{\Delta}{\sqrt{m_0^3\alpha_-}}\right)t\ . \end{eqnarray*}
Inequalities $0<\varepsilon_0<1$ and $t\le \min\left\{\sqrt{\frac{\alpha_-^3}{m_0}}\ ,\ \frac{\sqrt{m_0^3\alpha_-}}{\Delta} \right\}$ imply
\begin{eqnarray*} t\le\frac{2}{2\max\left\{\sqrt{\frac{m_0}{\alpha_-^3}},\ \frac{\Delta}{\sqrt{m_0^3\alpha_-}}\right\}}
\le\frac{\pi-\sqrt{\varepsilon_0}}{\sqrt{\frac{m_0}{\alpha_-^3}}+\frac{\Delta}{\sqrt{m_0^3\alpha_-}}}\quad \Longrightarrow\quad |x_*(t)-x_*(0)|\le \pi-\sqrt{\varepsilon_0}
\end{eqnarray*}
 and hence $x_*(t)\in {\mathbb X}$ for $0\le t\le T$. 
Since $0\le t\le 2^N\frac{s_0\delta}{\Delta}$, we have \begin{eqnarray*}|{\cal G}_*(t)-{\cal G}_*(0)|\le \frac{2^{-(N+1)}\Delta t}{s_0}\le \frac{\delta}{2}\end{eqnarray*} and hence ${\cal G}_*(t)\in {\mathbb G}$.  Let us now evaluate the variation of $\gamma_*$ 
during the time $T$. We have
\begin{eqnarray*}
|\gamma_*(T)-\gamma_*(0)|&\ge& \inf|\partial_{{\cal G}_*}(g_*+f_*)|t\ge\left(\inf|\partial_{{\cal G}_*}\overline f|-\sup|\partial_{{\cal G}_*}(|g_*-\overline f|+|f_*|)| \right)T\nonumber\\
&\ge& \left(\inf|\partial_{{\cal G}_*}\overline f|-\frac{\Delta}{\Lambda} N_0^{-1} \right)T
\end{eqnarray*}
 {Proceeding as in~(\ref{Delta2}) and using Cauchy inequalities, one sees that} $\inf|\partial_{{\cal G}_*}\overline f|\ge c^*\beta_*{m_0^2}\frac{a}{\Lambda\alpha_+^2}$. So, using~(\ref{Delta2}) and $N_0^{-1}<\frac{c_0^2\varepsilon_0^2\alpha_-^2}{2\alpha_+^2}$,
\begin{eqnarray*}
|\gamma_*(T)-\gamma_*(0)|&\ge& 
c^*\beta_*{m_0^2}\frac{a}{\Lambda\alpha_+^2}\left(1-\frac{\alpha_+^2}{c_0^2\varepsilon_0^2\alpha_-^2N_0}\right)T\ge\frac{c^*}{2}\beta_*{m_0^2}\frac{a}{\Lambda\alpha_+^2}T\nonumber\\
&=&\frac{c^*}{2}m_0^2\beta_*\frac{a}{\Lambda\alpha_+^2}\min\left\{\sqrt{\frac{\alpha_-^3}{m_0}}\ ,\   \frac{\sqrt{m_0^3\alpha_-\varepsilon_0}}{\Delta} \ ,\ 2^N\frac{s_0\delta}{\Delta}\right\}\nonumber\\
&=&c^\circ \min\left\{\beta_*\sqrt{\frac{\alpha_-^3 a}{\alpha_+^4}}\ ,\ \frac{c^2_0\varepsilon^{5/2}_0\alpha^2_-}{\alpha_+^2}\sqrt{\frac{\alpha_-}{a}}\ ,\ 
\frac{c^2_0\varepsilon^{2}_0\alpha^2_-}{\alpha_+^2}2^{N_0} s_0\frac{\delta}{\Lambda}
\right\}{=:\frac{3\pi}{\eta}}
\end{eqnarray*}
with $c^*$, $c^\circ$ independent of $\alpha_-$, $\alpha_+$, $\beta$, $\overline\beta$, $\delta$, $\varepsilon_0$ and $s_0$.
We then see that   $|\gamma_*(T)-\gamma_*(0)|$ is lower bounded by $3\pi$ as soon as the  condition in~(\ref{last cond}) is satisfied. $\qquad \square$

 \appendix
\section{Proof of Theorem~\ref{NFL}}\label{Time--one flows}
In this section we  use  the definitions and the notations introduced in Section~\ref{A normal form lemma without small divisors}, plus the following further ones.

\begin{definition}\rm
Given ${\rm h}$ and $f$ as in~(\ref{h0f0}), we call\footnote{{\sc nqp} stands for ``non quasi--periodic''.} the function~(\ref{phiNEW}) {\it {\sc nqp}--primitive of $f$ relatively to ${\rm h}$% and ${\mathbb L}$
} or, simply, \it{{\sc nqp}--primitive of $f$}.

\end{definition}

\begin{example}\rm
Let $n=1$ and $f(I, \varphi)=a(I)\cos\varphi$. The {\sc nqp}--primitive of $f$ is \begin{eqnarray*}
\phi(I, \varphi, y, x)&=&\frac{ 1}{\omega_y}\int_0^x a(I)\cos\left(\varphi+\frac{\omega_{ I}}{\omega_y}(\tau-x)\right) d\tau=\frac{ a(I)}{\omega_I}\int_{-\frac{\omega_{ I}}{\omega_y}x}^0 \cos(\varphi+\psi)d\psi\nonumber\\
&=&\frac{a(I)}{\omega_I}\left(
\sin\varphi-\sin\left(\varphi-\frac{\omega_{ I}}{\omega_y}x\right)
\right)
\end{eqnarray*}
having changed variable $\tau\to\psi=\frac{\omega_{ I}}{\omega_y}(\tau-x) $.
\end{example}

\begin{definition}[time--one  flow]\rm
Let ${\cal L}_\phi(\cdot):=\big\{\phi, \cdot\big\}$. For a given $\phi\in {\cal O}_{\rho, s, \delta, r, \xi}$ and $h\in {\mathbb N}$, we denote as $\Phi_h$  the formal series

\begin{eqnarray} \label{queue}\Phi_h:=\sum_{j\ge h}\frac{{\cal L}_\phi^j}{j!}\ .\end{eqnarray} 
We call $\Psi:=\Phi_0:=e^{{\cal L}_\phi}$ {\it time--one flow} generated by $\phi$.
\end{definition}

\begin{definition}
\rm We call {\it {\sc nqp}--homological transformation} the time--one flow
$\Psi:=\Phi_0=e^{{\cal L}_\phi}$
generated by $\phi$ in~(\ref{exp}), with $\phi_{khj}$ as in~(\ref{solution}).
\end{definition}
Below we prove that the {\sc nqp}--homological transformation  is well defined. Later, we shall prove that  the composition of $N+1$ {\sc nqp}--homological transformations $\Psi_1=e^{{\cal L}_{\phi_1}}$, $\cdots$, $\Psi_{N+1}=e^{{\cal L}_{\phi_{N+1}}}$, where
\begin{eqnarray*}\Psi_j:\quad {\rm H}_{j-1}\to {\rm H}_{j} \end{eqnarray*}
with ${\rm H}_0={\rm H}$ provides the transformation $\Psi_*$ in~(\ref{PsiN}).

\begin{lemma}[\cite{poschel93}]\label{base lemma}
There exists an integer number $\overline{\rm c}_{n,m}$ such that, for any $\phi\in {\cal O}_{\rho, s, \delta, r, \xi}$ and any $r'<r$, $s'<s$, $\rho'<\rho$, $\xi'<\xi$, $\delta'<\delta$ such that
\begin{eqnarray*}\frac{\overline{\rm c}_{n,m}\|\phi\|_{\rho, s, \delta, r, \xi}}{d}<1\qquad d:=\min\big\{\rho'\sigma', r'\xi', {\delta'}^2\big\} \end{eqnarray*}
then  the series in~(\ref{queue}) converge uniformly so as to define 
the family $\{\Phi_h\}_{h=0,1,\cdots}$ of operators 
 \begin{eqnarray*}\Phi_h:\quad  {\cal O}_{\rho, s, \delta, r, \xi}\to {\cal O}_{\rho-\rho', s-s', \delta-\delta', r-r', \xi-\xi'}\ . \end{eqnarray*}
 Moreover, the following bound holds (showing, in particular, uniform  convergence):

\begin{eqnarray} \label{geometric series}\|{\cal L}^j_\phi[g]\|_{\rho-\rho', s-s', \delta-\delta', r-r', \xi-\xi'}\le j!\big(\frac{\overline{\rm c}_{n,m}\|\phi\|_{\rho, s, \delta, r, \xi}}{d}\big)^j\|g\|_{\rho, s, \delta, r, \xi}\end{eqnarray} 
for all $g\in{\cal O}_{\rho, s, \delta, r, \xi}$.
\end{lemma}

\begin{remark}[\cite{poschel93}]\rm The bound~(\ref{geometric series}) immediately implies
\begin{eqnarray} \label{h power}\|\Phi_hg\|_{\rho-\rho', s-s', \delta-\delta', r-r', \xi-\xi'}\le \frac{\big(\frac{\overline{\rm c}_{n,m}\|\phi\|_{\rho, s, \delta, r, \xi}}{d}\big)^h}{1-\frac{\overline{\rm c}_{n,m}\|\phi\|_{\rho, s, \delta, r, \xi}}{d}}\|g\|_{\rho, s, \delta, r, \xi}\qquad \forall g\in {\cal O}_{\rho, s, \delta, r, \xi}\ .\end{eqnarray} 
\end{remark}

\begin{lemma}[Iterative Lemma]\label{iterative lemma}
There exists a number $\widetilde{\rm c}_{n,m}>1$ such that the following holds. For any choice  of positive numbers  $r'$, $\rho'$, $s'$, $\xi'$, $\delta'$ satisfying
%\begin{eqnarray*}
%&&r'<r\ ,\quad \r'<\r\ ,\quad \xi'<\xi\\
%&& s'<s\ ,\quad \delta'<\d\ ,\quad {\cal X}\left\|\frac{\omega_{\mathbf I}}{\omega_y}\right\|_{\rho, r}<s-s'\ ,\quad
%{\cal X}\left\|\frac{\omega_{\mathbf J}}{\omega_y}\right\|_{\rho, r}<
%\log\frac{\delta}{\delta'} 
%\end{eqnarray*}
\begin{eqnarray} \label{ineq1}
&&{2\rho'<\rho\ ,\quad 2r'<r\ ,\quad  2\xi'<\xi}\\
\label{ineq2}
&&2s'<s\ ,\quad 2\delta'<\delta\ ,\quad { {\cal X}\left\|\Im\frac{\omega_{\mathbf I}}{\omega_y}\right\|_{\rho, r}<s-2s'\ ,\quad
{\cal X}\left\|\frac{\omega_{\mathbf J}}{\omega_y}\right\|_{\rho, r}<
\log\frac{\delta}{2\delta'} }
\end{eqnarray} 
and and provided that the  following inequality holds true
\begin{eqnarray} \label{smallness}
 \widetilde{\rm c}_{n,m}\frac{{\cal X}}{d}   % \|\frac{1}{\omega_y}\|_{\rho, r}
 \left \|\frac{\widetilde f}{\omega_y}\right\|_{\rho, s, \delta, r, \xi} <1\qquad d:=\min\big\{\rho'\sigma', r'\xi', {\delta'}^2\big\}
\end{eqnarray} 
 one can find an operator \begin{eqnarray*}\Phi:\quad {\cal O}_{\rho, s, \delta, r, \xi}\to {\cal O}_{\rho_+, s_+, \delta_+, r_+, \xi_+} \end{eqnarray*}
 with
%  $$r_+:=r-r'\ ,\quad \rho_+:=\rho-\rho'\ ,\quad \xi_+:=\xi-\xi'\ ,\quad s_+:=s-s'-{\cal X}\left\|\frac{\omega_{\mathbf I}}{\omega_y}\right\|_{\rho, r}\ ,\quad \d_+:=\delta e^{-{\cal X}\left\|\frac{\omega_{\mathbf J}}{\omega_y}\right\|_{\rho, r}}-\delta'$$
\begin{eqnarray*}
&&r_+:=r-2r'\ ,\quad \rho_+:=\rho-2\rho'\ ,\quad \xi_+:=\xi-2\xi'\ ,\quad s_+:=s-2s'-{\cal X}\left\|\Im\frac{\omega_{\mathbf I}}{\omega_y}\right\|_{\rho, r}\nonumber\\
&& \delta_+:=\delta e^{-{\cal X}\left\|\frac{\omega_{\mathbf J}}{\omega_y}\right\|_{\rho, r}}-2\delta'
\end{eqnarray*}
which carries the Hamiltonian  ${\rm H}$ in~(\ref{step i}) to
\begin{eqnarray*}{\rm H}_+:=\Phi[{\rm H}]={\rm h}+g+\overline f+f_+ \end{eqnarray*}
where 
\begin{eqnarray} \label{bound}\|f_+\|_{r_+,\rho_+, \xi_+, s_+,\delta_+}\le \widetilde{\rm c}_{n,m}\left(\frac{{\cal X}}{d}\left \|\frac{\widetilde f}{\omega_y}\right\|_{\rho, s, \delta, r, \xi}\| f\|_{\rho, s, \delta, r, \xi}+ \|\{\phi, g\}\|_{\rho_1-\rho', s_1-s', \delta_1-\delta', r_1-r', \xi_1-\xi'}\right)
\end{eqnarray} 
with \begin{eqnarray*}\rho_1:=\rho\ ,\quad s_1:=s-{\cal X}\left\|\Im\frac{\omega_{\mathbf I}}{\omega_y}\right\|_{\rho, r}\ ,\quad \delta_1:=\delta e^{-{\cal X}\left\|\frac{\omega_{\mathbf J}}{\omega_y}\right\|_{\rho, r}}\ ,\quad r_1:=r\ ,\quad  \xi_1:=\xi \end{eqnarray*} for a suitable $\phi\in {\cal O}_{\rho_1, s_1, \delta_1, r_1, \xi_1}$ verifying \begin{eqnarray} \label{bound on phi}\|\phi\|_{\rho_1, s_1, \delta_1, r_1, \xi_1 }\le {\cal X}   %\|\frac{1}{\omega_y}\|_{\rho, r} \|\widetilde f\|_{\rho, s, \delta, r, \xi}
\left \|\frac{\widetilde f}{\omega_y}\right\|_{\rho, s, \delta, r, \xi}
\ .\end{eqnarray} 
Furthermore, 
if
\begin{eqnarray*}({\mathbf I}_+, \bm\varphi_+, {\mathbf p}_+, {\mathbf q}_+, y_+, x_+):=\Phi({\mathbf I}, \bm\varphi, {\mathbf p}, {\mathbf q}, y, x) \end{eqnarray*}
the following uniform bounds hold:
\begin{eqnarray} \label{phi close to id}
&&\max\Big\{s'|{\mathbf I}-{\mathbf I}_+|,\ \rho'|\bm\varphi -\bm\varphi _+|,\ 
\delta' |{\mathbf p}-{\mathbf p}_+|,\ \delta |{\mathbf q}-{\mathbf q}_+|,\ \xi' |y-y_+|,\ r' |x-x_+| \Big\}\nonumber\\
&&\leq2\, {\cal X}   % \|\frac{1}{\omega_y}\|_{\rho, r} \|\widetilde f\|_{\rho, s, \delta, r, \xi}
\left \|\frac{\widetilde f}{\omega_y}\right\|_{\rho, s, \delta, r, \xi}
\ .
\end{eqnarray} 
\end{lemma}

\begin{remark}
\rm 
The right hand side of~(\ref{bound}) benefits of the absence of small divisors, as no ``ultraviolet'' (i.e,  with size $\sim e^{-Ns}\|f\|$) term  appears.
\end{remark}

\par\medskip\noindent{\bf Proof\ } 
Let $\overline{\rm c}_{n,m}$ be as in Lemma~\ref{base lemma}. We shall choose $\widetilde{\rm c}_{n,m}$ suitably large with respect to $\overline{\rm c}_{n,m}$.

{\smallskip\noindent}
Let $\phi_{khj}$ as in~(\ref{solution}). 
Let us fix
\begin{eqnarray} \label{ovl param}0<\overline r\le r\ ,\quad 0<\overline\rho\le \rho\ ,\quad 0<\overline\xi\le \xi\ ,\quad 0<\overline s< s\ ,\quad 0<\overline\delta< \delta\end{eqnarray} 
and assume that
\begin{eqnarray} \label{ovl param1}{\cal X}\left\|\Im\frac{\omega_{\mathbf I}}{\omega_y}\right\|_{\rho, r} \le s-\overline s\ ,\qquad {\cal X}\left\|\frac{\omega_{\mathbf J}}{\omega_y}\right\|_{\rho, r} \le \log\frac{\delta}{\overline \delta}\ .\end{eqnarray} 
Then we have
\begin{eqnarray*}\|\phi_{khj}\|_{\overline\rho,\overline r,\overline\xi}\le  %\|\frac{1}{\omega_y}\|_{\overline\rho, \overline r}
 \left\|\frac{f_{khj}}{\omega_y}\right\|_{\overline\rho,\overline r,\overline\xi}\left\|\int_0^x |e^{-\frac{\l_{khj}}{\omega_y}\tau}|\right\|_{\overline\rho,\overline r,\overline\xi}d\tau\le{\cal X}   %\|\frac{1}{\omega_y}\|_{\overline\rho, \overline r} 
  \left\|\frac{f_{khj}}{\omega_y}\right\|_{\overline\rho,\overline r,\overline\xi} e^{{\cal X}\left\|\frac{\l_{khj}}{\omega_y}\right\|_{\overline\rho, \overline r}}\ . \end{eqnarray*}
%We take, in particular, $\overline r$, $\overline\rho$, ${\cal X}$ as above and, moreover,
%$\overline s$, $\overline\delta$ such that
% $$\overline s+{\cal X}\left\|\frac{\omega_{\mathbf I}}{\omega_y}\right\|_{\overline\rho, \overline r} \le s\ ,\qquad \overline \delta e^{{\cal X}\left\|\frac{\omega_{\mathbf J}}{\omega_y}\right\|_{\overline\rho, \overline r} }\le \d\ .$$
Since
\begin{eqnarray*}\|\frac{\l_{khj}}{\omega_y}\|_{\overline\rho, \overline r} \le (h+j)\left\|\frac{\omega_{\mathbf J}}{\omega_y}\right\|_{\overline\rho, \overline r} +|k|\left\|\Im\frac{\omega_{\mathbf I}}{\omega_y}\right\|_{\overline\rho, \overline r} 
\end{eqnarray*}
we have
\begin{eqnarray*}\|\phi_{khj}\|_{\overline\rho,\overline r,\overline\xi}\le {\cal X}   %\|\frac{1}{\omega_y}\|_{\overline\rho, \overline r}  
\left\|\frac{\widetilde f_{khj}}{\omega_y}\right\|_{\overline\rho,\overline r,\overline\xi}
e^{(h+j){\cal X}\left\|\frac{\omega_{\mathbf J}}{\omega_y}\right\|_{\overline\rho, \overline r} +|k|{\cal X}\left\|\Im\frac{\omega_{\mathbf I}}{\omega_y}\right\|_{\overline\rho, \overline r} }\ .
\end{eqnarray*}
which yields (after multiplying by $e^{|k|\overline s}(\overline\delta)^{j+h}$ and summing over $k$, $j$, $h$ with $(k,h-k)\ne (0,0)$) to
\begin{eqnarray*}\|\phi\|_{\overline\rho,\overline r,\overline\xi, \overline s,\overline\delta}\le{\cal X}  %\|\frac{1}{\omega_y}\|_{\overline r,\overline \r} 
\left\|\frac{\widetilde f}{\omega_y}\right\|_{\overline r,\overline\rho,\overline \xi,\overline s+{\cal X}\left\|\Im\frac{\omega_{\mathbf I}}{\omega_y}\right\|_{\overline\rho, \overline r} ,\overline \delta e^{{\cal X}\left\|\frac{\omega_{\mathbf J}}{\omega_y}\right\|_{\overline\rho, \overline r} }}\ . \end{eqnarray*}
Note that the right hand side is well defined because of~(\ref{ovl param1}).
In the case of the choice
\begin{eqnarray*}
&&\overline r=r=:r_1\ ,\quad\overline\rho=\rho=: \rho_1\ ,\quad \overline\xi=\xi=:\xi_1\ ,\quad \overline s=s-{\cal X}\left\|\Im\frac{\omega_{\mathbf I}}{\omega_y}\right\|_{\rho, r}=:s_1\nonumber\\
&& \overline\delta=\delta e^{-{\cal X}\left\|\frac{\omega_{\mathbf J}}{\omega_y}\right\|_{\rho, r}}=:\delta_1
\end{eqnarray*}
(which, in view of the two latter inequalities in~(\ref{ineq2}), satisfies~(\ref{ovl param})--(\ref{ovl param1})) the inequality becomes~(\ref{bound on phi}).
An application of Lemma~\ref{base lemma},with $r$, $\rho$, $\xi$, $s$, $\delta$ replaced by $r_1-r'$, $\rho_1-\rho'$, $\xi_1-\xi'$, $s_1-s'$, $\delta_1-\delta'$,  concludes  with a suitable choice of $\widetilde{\rm c}_{n,m}>\overline{\rm c}_{n,m}$ and (by~(\ref{f1})) \begin{eqnarray*} f_+:=\Phi_2({\rm h})+\Phi_1(g)+\Phi_1(f)\ . \end{eqnarray*}
Observe that the bound~(\ref{bound}) follows from  Equations~(\ref{h power}),~(\ref{geometric series}) and the identities
\begin{eqnarray*}\Phi_2[{\rm h}]=\sum_{j=2}^\infty \frac{{\cal L}^j_\phi({\rm h})}{j!}=\sum_{j=1}^\infty \frac{{\cal L}^{j+1}_\phi({\rm h})}{(j+1)!}=-\sum_{j=1}^\infty \frac{{\cal L}^{j}_\phi(\widetilde f)}{(j+1)!}
\end{eqnarray*}
\begin{eqnarray*}\Phi_1[g]=\sum_{j=1}^\infty \frac{{\cal L}^j_\phi(g)}{j!}=\sum_{j=0}^\infty \frac{{\cal L}^{j+1}_\phi(g)}{(j+1)!}=\sum_{j=0}^\infty \frac{{\cal L}^{j}_\phi(g_1)}{(j+1)!} \end{eqnarray*}
with $g_1:={\cal L}_\phi(g)=\{\phi, g\}$.
The bounds in~(\ref{phi close to id}) are a consequence of equalities of the kind
\begin{eqnarray*}{\mathbf I}_+-{\mathbf I}=\sum_{j=0}^\infty\frac{{\cal L}^{j+1}_\phi({\mathbf I})}{(j+1)!}=\sum_{j=0}^\infty\frac{{\cal L}^{j}_\phi(-\partial_\varphi\phi)}{(j+1)!} \end{eqnarray*}
(and similar).
 $\qquad\square$

{\smallskip\noindent}
The proof of the  Theorem~\ref{NFL} goes through iterate applications of Lemma~\ref{iterative lemma}.\\ We premise the following

\begin{remark}\label{stronger iterative lemma}\rm
Replacing conditions in~(\ref{ineq2}) with the stronger ones\footnote{The three first inequalities in~(\ref{new cond}) are immediately seen to be stronger that the corresponding three first inequalities  in~(\ref{ineq2}).
On the other hand, rewriting the second inequality in~(\ref{new cond}) as
$\frac{\delta'}{\delta}<1-\frac{2\delta'}{\delta}$
and using the inequality (which holds for all $x\ge 1$) $\log x\ge 1-\frac{1}{x}$  with $x=\frac{\delta}{2\delta'}$, we have  also $\frac{\delta'}{\delta}<\log\frac{\delta}{2\delta'} $. 
}
%\beqno2s'<s\ ,\quad 2\delta'<\d\ ,\quad {\cal X}\left\|\frac{\omega_{\mathbf I}}{\omega_y}\right\|_{\rho, r}<s'\ ,\quad
%{\cal X}\left\|\frac{\omega_{\mathbf J}}{\omega_y}\right\|_{\rho, r}<
%\frac{\delta'}{\delta} {\rm e}qno

\begin{eqnarray} \label{new cond}
{3s'<s\ ,\quad 3\delta'<\delta\ ,\quad {\cal X}\left\|\Im\frac{\omega_{\mathbf I}}{\omega_y}\right\|_{\rho, r}<s'\ ,\quad
{\cal X}\left\|\frac{\omega_{\mathbf J}}{\omega_y}\right\|_{\rho, r}<
\frac{\delta'}{\delta}} \end{eqnarray} 

{\smallskip\noindent}
(and keeping~(\ref{ineq1}),~(\ref{smallness}) unvaried)  one can take, for $s_+$, $\delta_+$,  $s_1$, $\delta_1$ the simpler expressions
%$$s_{+\rm new}=s-2s'\ ,\quad \d_{+\rm new}=\d-2\delta'\ ,\quad s_{1\rm new}:=s-s'\ ,\quad \d_{1\rm new}=\d-\delta'\ $$

\begin{eqnarray*}{s_{+\rm new}=s-3s'\ ,\quad \delta_{+\rm new}=\delta-3\delta'\ ,\quad s_{1\rm new}:=s-s'\ ,\quad \delta_{1\rm new}=\delta-\delta'\ } \end{eqnarray*}
(while keeping $r_+$, $\rho_+$, $\xi_+$, $r_1$, $\rho_1$, $\xi_1$ unvaried).
Indeed, since $1-e^{-x}\le x$ for all $x$,
\begin{eqnarray*}\delta_1=\delta e^{-{\cal X}\left\|\frac{\omega_{\mathbf J}}{\omega_y}\right\|_{\rho, r}} = \delta-\delta(1- e^{-{\cal X}\left\|\frac{\omega_{\mathbf J}}{\omega_y}\right\|_{\rho, r}})\ge \delta-\delta{\cal X}\left\|\frac{\omega_{\mathbf J}}{\omega_y}\right\|_{\rho, r}\ge \delta-\delta' =\delta_{1\rm new}\ . \end{eqnarray*}
This also implies $\xi_+=\delta_1-\delta'\ge \delta-2\delta'=\xi_{+\rm new}$. That $s_+\ge s_{+\rm new}$, $s_1\ge s_{1\rm new}$ is even more immediate.\end{remark}

{\smallskip\noindent}
Now we can proceed with the

\paragraph{Proof of    Theorem~\ref{NFL}}
Let $\widetilde{\rm c}_{n,m}$ be as in Lemma~\ref{iterative lemma}. We shall choose $\overline{\rm c}_{n,m}$ suitably large with respect to $\widetilde{\rm c}_{n,m}$.

{\smallskip\noindent}
We apply Lemma~\ref{iterative lemma} with
\begin{eqnarray*}{2}\rho'=\frac{\rho}{3}\ ,\quad {3}s'=\frac{s}{3}\ ,\quad {3}\delta'=\frac{\delta}{3}\ ,\quad{2}r'=\frac{r}{3}\ ,\quad {2}\xi'=\frac{\xi}{3}\ ,\quad g\equiv0\ . \end{eqnarray*}
We make use of the stronger formulation described in Remark~\ref{stronger iterative lemma}.  Conditions in~(\ref{ineq1}) and the two former conditions in~(\ref{new cond}) are trivially true. The two latter inequalities in 
(\ref{new cond}) reduce to 
\begin{eqnarray*}{\cal X}\left\|\Im\frac{\omega_{\mathbf I}}{\omega_y}\right\|_{\rho, r}<\frac{s}{9}\ ,\quad
{\cal X}\left\|\frac{\omega_{\mathbf J}}{\omega_y}\right\|_{\rho, r}<
\frac{1}{9} \end{eqnarray*}
and they are certainly satisfied by assumption~(\ref{normal form assumptions}), for {$N>2$}. Since  
 \begin{eqnarray*} d=\min\{ \rho' s', r'\xi', {\delta'}^2 \}=\min\{ \rho s/{36}, r\xi/{54}, {\delta}^2/{81} \}\ge\frac{1}{{81}}
\min\{ \rho s, r\xi, {\delta}^2 \}=\frac{\rm d}{{81}}
\end{eqnarray*}
we have that condition~(\ref{smallness}) is certainly implied by the last inequality in~(\ref{normal form assumptions}), once one chooses ${\rm c}_{n,m}>{81} \widetilde{\rm c}_{n,m}$.
By Lemma~\ref{iterative lemma}, it is then possible to conjugate ${\rm H}$ to
\begin{eqnarray*}{\rm H}_1={\rm h}_0+\overline f+f_1\end{eqnarray*}
with $f_1\in {\cal O}_{\rho^{(1)},s^{(1)},\delta^{(1)}, r^{(1)},\xi^{(1)},}$, where $(\rho^{(1)},s^{(1)},\delta^{(1)}, r^{(1)},\xi^{(1)},):=2/3 (\rho, s, \delta, r, \xi)$ and
\begin{eqnarray} \label{f1***}\|f_1\|_{\rho^{(1)},s^{(1)},\delta^{(1)}, r^{(1)},\xi^{(1)},}\le {81}\widetilde{\rm c}_{n,m} \frac{{\cal X}}{\rm d}  % \|\frac{1}{\omega_y}\|_{\rho, r} 
\left\|\frac{\widetilde f}{\omega_y}\right\|_{\rho, s, \delta, r, \xi}\| f\|_{\rho, s, \delta, r, \xi}\le \frac{\| f\|_{\rho, s, \delta, r, \xi}}{2}\end{eqnarray} 
since ${\rm c}_{n,m}\ge {162} \widetilde{\rm c}_{n,m}$ and $N\ge1$. Now we aim to apply Lemma~\ref{iterative lemma} $N$ times, again as described  in Remark~\ref{stronger iterative lemma}, each time with parameters
\begin{eqnarray} \label{values0} \rho_j'=\frac{\rho}{6N}\ ,\quad s_j'=\frac{s}{{9}N}\ ,\quad \delta_j'=\frac{\delta}{{9}N}\ ,\quad r_j'=\frac{r}{6N}\ ,\quad \xi_j'=\frac{\xi}{6N}\ .\end{eqnarray} 
Therefore we find, at each step
\begin{eqnarray} \label{values}
&&\rho^{(j+1)}:=\rho^{(1)}-j\frac{\rho}{3N}\ ,\quad s^{(j+1)}:=s^{(1)}-j\frac{s}{3N}\ ,\quad \delta^{(j+1)}:=\delta^{(1)}-j\frac{\delta}{3N}\nonumber\\
&&r^{(j+1)}:=r^{(1)}-j\frac{r}{3N}\ ,\quad \xi^{(j+1)}:=\xi^{(1)}-j\frac{\xi}{3N}\nonumber\\
&&\rho_1^{(j)}:=\rho^{(j)}\ ,\quad s_1^{(j)}:=s^{(j)}-\frac{s}{9N}\ ,\quad \delta_1^{(j)}:=\delta^{(j)}-\frac{\delta}{9N}\nonumber\\
&&r_1^{(j)}:=r^{(j)}\ ,\quad  \xi_1^{(j)}:=\xi^{(j)}\ ,\quad {\cal X}_j:=\sup\{|x|:\ x\in {\mathbb X}_{\xi_j}\}
\end{eqnarray} 
with $1\le j\le N$.

{\smallskip\noindent}
We assume that for a certain $1\le i\le N$ and all $1\le j\le i$, we have ${\rm H}_j\in {\cal O}_{\rho^{(j)},s^{(j)},\delta^{(j)}, r^{(j)},\xi^{(j)},}$ of the form
\begin{eqnarray} \label{Hi}
&&{\rm H}_j={\rm h}_0+g_{j-1}+f_j\ , \quad g_{j-1}\in {\cal N}_{\rho^{(j)},s^{(j)},\delta^{(j)}, r^{(j)},\xi^{(j)},}\ ,\quad g_{j-1}-g_{j-2}=\overline f_{j-1}\\
&&\label{i+1 ineq} \|f_j\|_{\rho^{(j)},s^{(j)},\delta^{(j)}, r^{(j)},\xi^{(j)},}\le \frac{\|f_1\|_{\rho^{(1)},s^{(1)},\delta^{(1)}, r^{(1)},\xi^{(1)},}}{2^{j-1}}\end{eqnarray} 
with  $g_{-1}\equiv0$,  $g_0=\overline f$. We  want to prove that, if $i<N$,
Lemma~\ref{iterative lemma} can be applied once again, so as to conjugate ${\rm H}_i$ to a suitable ${\rm H}_{i+1}$ such that~(\ref{Hi})--(\ref{i+1 ineq})
are true with $j=i+1$.
To this end, 
according to the discussion in Remark~\ref{stronger iterative lemma},  
we check the  stronger inequalities
\begin{eqnarray} \label{ineq0}
{2{\rho'_i}<\rho^{(i)}\ ,\quad 2r'_i<r^{(i)}\ ,\quad  2\xi'_i<\xi^{(i)}}\end{eqnarray} 
\begin{eqnarray} \label{smallness0}
&&{\cal X}_i\left\|\Im\frac{\omega_{\mathbf I}}{\omega_y}\right\|_{\rho_i, r_i}<s'_i%=\frac{s}{6N}
\ ,\quad
{\cal X}_i\left\|\frac{\omega_{\mathbf J}}{\omega_y}\right\|_{\rho_i, r_i}<
\frac{\delta'_i}{\delta_i}\\
\label{smallness i}
&& \widetilde{\rm c}_{n,m}\frac{{\cal X}_i}{d_i} % \|\frac{1}{\omega_y}\|_{\rho_i, r_i}
\left\|\frac{f_i}{\omega_y}\right\|_{\rho_i,s_i,\delta_i, r_i,\xi_i}<1\ .
\end{eqnarray} 
where
$d_i:=\min\{ \rho_i' s_i', r_i'\xi_i', {\delta'}_i^2 \}$.
Conditions~(\ref{ineq0}) and~(\ref{smallness0}) are certainly  verified, since in fact they are implied by the definitions above 
(using also $\delta_i\le \frac{2}{3}\delta$, ${\cal X}_i\le {\cal X}$) and the two former inequalities in~(\ref{normal form assumptions}).
 To check the validity of~(\ref{smallness i}), we firstly observe that
\begin{eqnarray} \label{di}d_i=\min\{\rho'^{(i)} s'^{(i)},\ (\delta'^{(i)})^2,\ r'^{(i)}\xi'_j,\}\ge \frac{\rm d}{{81}N^2}\ .\end{eqnarray} 
Using then ${\rm c}_{n,m}>{162} \widetilde{\rm c}_{n,m}$,${\cal X}_i<{\cal X}$, Equation~(\ref{f1***}), the inequality in~(\ref{i+1 ineq}) with $j=i$ and the last inequality in~(\ref{normal form assumptions}), we easily conclude
\begin{eqnarray*}
&&\|f_i\|_{\rho^{(i)},s^{(i)},\delta^{(i)}, r^{(i)},\xi^{(i)}}\le \|f_1\|_{\rho^{(1)},s^{(1)},\delta^{(1)}, r^{(1)},\xi^{(1)},}\le {81}\widetilde{\rm c}_{n,m} \frac{{\cal X}}{\rm d}   %\|\frac{1}{\omega_y}\|_{\rho, r} 
\left\|\frac{\widetilde f}{\omega_y}\right\|_{\rho, s, \delta, r, \xi}\|{ f}\|_{\rho, s, \delta, r, \xi}
 \nonumber\\
 &&\le
 \frac{1}{\widetilde{\rm c}_{n,m}}\frac{\rm d}{{81} N^2}\frac{1} {{\cal X}}  \left\|\frac{1}{\omega_y}\right\|_{\rho, r}^{-1}\le  \frac{1}{\widetilde{\rm c}_{n,m}}\frac{d_i} {{\cal X}_i}  \left\|\frac{1}{\omega_y}\right\|_{\rho^{(i)}, r^{(i)}}^{-1}\end{eqnarray*}
 which implies~(\ref{smallness i}).

{\smallskip\noindent}
 Then the Iterative Lemma is applicable to ${\rm H}_i$, which is conjugated to
 \begin{eqnarray} \label{fj+1}{\rm H}_{i+1}={\rm h}_0+g_{i}+f_{i+1}\ , \quad g_{i}\in {\cal N}_{\rho^{(i+1)},s^{(i+1)},\delta^{(i+1)}, r^{(i+1)},\xi^{(i+1)}}\ ,\quad g_{i}-g_{i-1}=\overline f_{i}\end{eqnarray} 
 with $g_{i}$, $f_{i+1}$ satisfying~(\ref{Hi})  with $j=i+1$.
We prove that~(\ref{i+1 ineq}) holds   with $j=i+1$, so as to complete the inductive step. By the thesis of the Iterative Lemma,
\begin{eqnarray*}
\left\|f_{i+1}\right\|_{\rho^{(i+1)}, s^{(i+1)}, \delta^{(i+1)}, r^{(i+1)}, \xi^{(i+1)}} &\le& \widetilde{\rm c}_{n,m}\Big(\frac{{\cal X}_i}{{d}_i}\left \|\frac{\widetilde f_i}{\omega_y}\right\|_{\rho^{(i)}, s^{(i)}, \delta^{(i)}, r^{(i)}, \xi^{(i)}}\| f_i\|_{\rho^{(i)}, s^{(i)}, \delta^{(i)}, r^{(i)}, \xi^{(i)}}\nonumber\\
&+& \|\{\phi_i, g_{i-1}\}\|_{\rho^{(i)}_1-{\rho'}^{(i)}, s^{(i)}_1-{s'}^{(i)}, \delta^{(i)}_1-{\delta'}^{(i)}, r^{(i)}_1-{r'}^{(i)}, \xi^{(i)}_1-{\xi'}^{(i)}}\Big)
\end{eqnarray*}
On the other hand, using~(\ref{f1***}),~(\ref{di}) and the last assumption in~(\ref{normal form assumptions}) and~(\ref{i+1 ineq})    with $j=i$, we obtain
\begin{eqnarray*}
\frac{{\cal X}_i}{{d}_i}\left \|\frac{\widetilde f_i}{\omega_y}\right\|_{\rho^{(i)}, s^{(i)}, \delta^{(i)}, r^{(i)}, \xi^{(i)}}&\le & \frac{81N^2\chi_0}{{\rm d}}\left\|\frac{1}{\omega_y}\right\|\left\|\widetilde f_i\right\|\le
\frac{81N^2\chi_0}{{\rm d}}\left\|\frac{1}{\omega_y}\right\|\left\| f_i\right\|
\nonumber\\
&\le & \frac{81N^2\chi_0}{{\rm d}}\left\|\frac{1}{\omega_y}\right\|\left\| f_1\right\|\nonumber\\
&\le & \frac{(81)^2 \widetilde{\rm c}_{n, m} N^2\chi^2_0}{{\rm d}^2}\left\|\frac{1}{\omega_y}\right\|\left\|\frac{\widetilde f}{\omega_y}\right\|_{\rho, s, \delta, r, \xi}\| f\|_{\rho, s, \delta, r, \xi}
\nonumber\\
&\le&\frac{1}{6\widetilde{\rm c}_{n, m}}
\end{eqnarray*}
Furthermore, using~(\ref{geometric series}) with $j=1$,~(\ref{bound on phi}) and
\begin{eqnarray*}
\|g_0\|_{\rho^{(i)}_1-{\rho'}^{(i)}, s^{(i)}_1-{s'}^{(i)}, \delta^{(i)}_1-{\delta'}^{(i)}, r^{(i)}_1-{r'}^{(i)}, \xi^{(i)}_1-{\xi'}^{(i)}}&\le& \|f\|_{\rho^{(i)}_1-{\rho'}^{(i)}, s^{(i)}_1-{s'}^{(i)}, \delta^{(i)}_1-{\delta'}^{(i)}, r^{(i)}_1-{r'}^{(i)}, \xi^{(i)}_1-{\xi'}^{(i)}} \nonumber\\
&\le&\|f\|_{\rho^{(0)}, s^{(0)}, \delta^{(0)}, r^{(0)}, \xi^{(0)}}
\end{eqnarray*}
and, with $j=1$, $\cdots$, $i-1$
\begin{eqnarray} \label{deltag}
\|g_j-g_{j-1}\|_{\rho^{(i)}_1-{\rho'}^{(i)}, s^{(i)}_1-{s'}^{(i)}, \delta^{(i)}_1-{\delta'}^{(i)}, r^{(i)}_1-{r'}^{(i)}, \xi^{(i)}_1-{\xi'}^{(i)}}&\le& \|f_j\|_{\rho^{(i)}_1-{\rho'}^{(i)}, s^{(i)}_1-{s'}^{(i)}, \delta^{(i)}_1-{\delta'}^{(i)}, r^{(i)}_1-{r'}^{(i)}, \xi^{(i)}_1-{\xi'}^{(i)}} \nonumber\\
&\le&\|f_j\|_{\rho^{(j)},s^{(j)},\delta^{(j)}, r^{(j)},\xi^{(j)},}\nonumber\\
&\le&\frac{\|f_1\|_{\rho^{(1)},s^{(1)},\delta^{(1)}, r^{(1)},\xi^{(1)},}}{2^{j-1}}
\end{eqnarray} 
we obtain
\begin{eqnarray*}
&&\|\{\phi_i, g_{i-1}\}\|_{\rho^{(i)}_1-{\rho'}^{(i)}, s^{(i)}_1-{s'}^{(i)}, \delta^{(i)}_1-{\delta'}^{(i)}, r^{(i)}_1-{r'}^{(i)}, \xi^{(i)}_1-{\xi'}^{(i)}}\nonumber\\
&&\qquad\le \|\{\phi_i, g_{0}\}\|_{\rho^{(i)}_1-{\rho'}^{(i)}, s^{(i)}_1-{s'}^{(i)}, \delta^{(i)}_1-{\delta'}^{(i)}, r^{(i)}_1-{r'}^{(i)}, \xi^{(i)}_1-{\xi'}^{(i)}}\nonumber\\
&&\qquad\qquad +\sum_{j=1}^{i-1}\|\{\phi_i, g_{j}-g_{j-1}\}\|_{\rho^{(i)}_1-{\rho'}^{(i)}, s^{(i)}_1-{s'}^{(i)}, \delta^{(i)}_1-{\delta'}^{(i)}, r^{(i)}_1-{r'}^{(i)}, \xi^{(i)}_1-{\xi'}^{(i)}}\nonumber\\
&&\qquad\qquad\le 2\frac{\overline{\rm c}_{n, m}\chi_i}{d_i}\left\|\frac{\widetilde f_i}{\omega_y}\right\|_{\rho^{(i)}, s^{(i)}, \delta^{(i)}, r^{(i)}, \xi^{(i)}}\Big(\frac{1}{N}
\|g_0\|_{\rho^{(i)}_1-{\rho'}^{(i)}, s^{(i)}_1-{s'}^{(i)}, \delta^{(i)}_1-{\delta'}^{(i)}, r^{(i)}_1-{r'}^{(i)}, \xi^{(i)}_1-{\xi'}^{(i)}}\nonumber\\
&&\qquad\qquad +\sum_{j=1}^{i-1}\|g_{j}-g_{j-1}\|_{\rho^{(i)}_1-{\rho'}^{(i)}, s^{(i)}_1-{s'}^{(i)}, \delta^{(i)}_1-{\delta'}^{(i)}, r^{(i)}_1-{r'}^{(i)}, \xi^{(i)}_1-{\xi'}^{(i)}}
\Big)\nonumber\\
&&\qquad\qquad\le 162 N^2\frac{\overline{\rm c}_{n, m}\chi_i}{{\rm d}}\left\|\frac{1}{\omega_y}\right\|%_{\rho^{(i)}, s^{(i)}, \delta^{(i)}, r^{(i)}, \xi^{(i)}}
\left(\frac{\|f\|_{\rho, s, \delta, r, \xi}}{N}+2\times {81}\widetilde{\rm c}_{n,m} \frac{{\cal X}}{\rm d}   %\|\frac{1}{\omega_y}\|_{\rho, r} 
\left\|\frac{\widetilde f}{\omega_y}\right\|_{\rho, s, \delta, r, \xi}\|{ f}\|_{\rho, s, \delta, r, \xi}
\right)\nonumber\\
&&\qquad\qquad\times
\left\|f_i\right\|_{\rho^{(i)}, s^{(i)}, \delta^{(i)}, r^{(i)}, \xi^{(i)}}\le \frac{1}3\left\|f_i\right\|_{\rho^{(i)}, s^{(i)}, \delta^{(i)}, r^{(i)}, \xi^{(i)}}
\end{eqnarray*}
where the last summand is to taken into account only when $i\ge 2$ and we have used that, when $j=0$, $d_i$ can be replaced with $N d_i$. Collecting the bounds above, we find~(\ref{i+1 ineq}) holds   with $j=i+1$.
Then the Iterative Lemma can be applied $N$ times, and we get~(\ref{thesis}) as a direct consequence of~(\ref{f1***}) and~(\ref{deltag}) and of the validity of~(\ref{i+1 ineq}) for all $1\le i\le N$. The inequalities in~(\ref{phi close to id***})
  follow from~(\ref{phi close to id}) and~(\ref{values0}). For example,
  \begin{eqnarray*}
  |{\mathbf I}_*-{\mathbf I}_0|\le\sum_{j=1}^N|{\mathbf I}^{(j)}-{\mathbf I}^{(j-1)}|\le2\frac{\chi_{0}}{{s'}^{(0)}}\left\|\frac{f_0}{\omega_y}\right\|+ 4\frac{\chi_{1}}{{s'}^{(1)}}\left\|\frac{1}{\omega_y}\right\|\|f_1\|\le 19\frac{\chi_{0}}{{s}^{(0)}}\left\|\frac{f_0}{\omega_y}\right\|
  \end{eqnarray*}
etc. $\qquad \square$

 \bibliographystyle{plain}
%\bibliography{REFERENCES.bib}
\addcontentsline{toc}{section}{References}
\def\cprime{$'$} \def\cprime{$'$}

\end{document}